\documentclass[a4paper,12pt]{amsart}
\usepackage{amssymb,amsthm}
\usepackage[all]{xy}

\begin{document}

\newcommand{\opp}{\bowtie }
\newcommand{\po}{\text {\rm pos}}
\newcommand{\supp}{\text {\rm supp}}
\newcommand{\End}{\text {\rm End}}
\newcommand{\diag}{\text {\rm diag}}
\newcommand{\Lie}{\text {\rm Lie}}
\newcommand{\Ad}{\text {\rm Ad}}
\newcommand{\car}{\mathcal R}
\newcommand{\Tr}{\rm Tr}
\newcommand{\Spec}{\text{\rm Spec}}

\def\ge{\geqslant}
\def\le{\leqslant}
\def\a{\alpha}
\def\b{\beta}
\def\c{\chi}
\def\g{\gamma}
\def\G{\Gamma}
\def\d{\delta}
\def\D{\Delta}
\def\L{\Lambda}
\def\e{\epsilon}
\def\et{\eta}
\def\io{\iota}
\def\o{\omega}
\def\p{\pi}
\def\ph{\phi}
\def\ps{\psi}
\def\r{\rho}
\def\s{\sigma}
\def\t{\tau}
\def\th{\theta}
\def\k{\kappa}
\def\l{\lambda}
\def\z{\zeta}
\def\v{\vartheta}
\def\va{\varphi}
\def\x{\xi}
\def\i{^{-1}}

\def\mapright#1{\smash{\mathop{\longrightarrow}\limits^{#1}}}
\def\mapleft#1{\smash{\mathop{\longleftarrow}\limits^{#1}}}
\def\mapdown#1{\Big\downarrow\rlap{$\vcenter{\text{$\scriptstyle#1$}}$}}
\def\mapup#1{\Big\uparrow\rlap{$\vcenter{\text{$\scriptstyle#1$}}$}}

\def\ca{\mathcal A}
\def\cb{\mathcal B}
\def\cc{\mathcal C}
\def\cd{\mathcal D}
\def\ce{\mathcal E}
\def\cf{\mathcal F}
\def\cg{\mathcal G}
\def\ch{\mathcal H}
\def\ci{\mathcal I}
\def\cj{\mathcal J}
\def\ck{\mathcal K}
\def\cl{\mathcal L}
\def\cm{\mathcal M}
\def\cn{\mathcal N}
\def\co{\mathcal O}
\def\cp{\mathcal P}
\def\cq{\mathcal Q}
\def\car{\mathcal R}
\def\cs{\mathcal S}
\def\ct{\mathcal T}
\def\cu{\mathcal U}
\def\cv{\mathcal V}
\def\cw{\mathcal W}
\def\cz{\mathcal Z}
\def\cx{\mathcal X}
\def\cy{\mathcal Y}

\def\tz{\tilde Z}
\def\tl{\tilde L}
\def\tc{\tilde C}
\def\ta{\tilde A}
\def\tb{\tilde B}
\def\tx{\tilde X}

\newtheorem*{th1}{Lemma 2.4}
\newtheorem*{th2}{Proposition 2.5}
\newtheorem*{th3}{Corollary 2.6}
\newtheorem*{th4}{Corollary 2.7}
\newtheorem*{th5}{Lemma 2.8}
\newtheorem*{th6}{Proposition 2.9}
\newtheorem*{th7}{Lemma 3.3}
\newtheorem*{th8}{Proposition 3.4}
\newtheorem*{th9}{Proposition 3.7}
\newtheorem*{th10}{Proposition 3.9}
\newtheorem*{th11}{Lemma 3.10}
\newtheorem*{th12}{Lemma 3.11}
\newtheorem*{th13}{Lemma 3.12}
\newtheorem*{th14}{Corollary 3.13}
\newtheorem*{th15}{Proposition 3.15}
\newtheorem*{th16}{Lemma 3.16}
\newtheorem*{th17}{Lemma 3.17}
\newtheorem*{th18}{Lemma 3.19}
\newtheorem*{th19}{Lemma 3.22}
\newtheorem*{th20}{Theorem 3.23}
\newtheorem*{th21}{Lemma 3.26}
\newtheorem*{th22}{Lemma 3.27}
\newtheorem*{th23}{Corollary 3.28}
\newtheorem*{th24}{Corollary 3.29}
\newtheorem*{th241}{Theorem 3.31}
\newtheorem*{th25}{Proposition 4.5}
\newtheorem*{th26}{Proposition 4.6}
\newtheorem*{th27}{Proposition 4.8}
\newtheorem*{th28}{Lemma 4.13}
\newtheorem*{th29}{Lemma 4.16}
\newtheorem*{th30}{Proposition 4.17}
\newtheorem*{th31}{Lemma 4.19}
\newtheorem*{th32}{Proposition 4.20}
\newtheorem*{th33}{Theorem 4.21}
\newtheorem*{th34}{Corollary 4.22}
\newtheorem*{th35}{Proposition 5.2}
\newtheorem*{th36}{Proposition 5.3}

\newtheorem{theorem}{Theorem}[section]
\newtheorem{lem}[theorem]{Lemma}
\newtheorem{cor}[theorem]{Corollary}
\newtheorem{prop}[theorem]{Proposition}
\newtheorem{thm}[theorem]{Theorem}
\newtheorem*{rmk}{Remark}
\newtheorem{eg}[theorem]{Example}

\title[Character sheaves on certain spherical varieties]
{Character sheaves on certain spherical varieties}
\author{Xuhua He}
\address{School of Mathematics, Institute for Advanced Study, Princeton, NJ 08540, USA}%
\thanks{The author is supported by NSF grant DMS-0111298.}
\email{hugo@math.ias.edu}


\begin{abstract}
We study a class of perverse sheaves on some spherical varieties
which include the strata of the De Concini-Procesi completion of a
symmetric variety. This is a generalization of the theory of
(parabolic) character sheaves.
\end{abstract}
\maketitle

\section*{Introduction}

\subsection*{0.1} Let $G$ be a connected, reductive algebraic group over an algebraically closed
field $\mathbf k$. In \cite{L3} and \cite{L4}, Lusztig introduced
the $(G \times G)$-varieties $Z_{J, y, D}$ and a class of
$G$-equivariant simple perverse sheaves on $Z_{J, y, D}$ which are
called ``(parabolic) character sheaves''. (The precise definition of
$Z_{J, y, D}$ can be found in 1.2 below). The varieties $Z_{J, y,
D}$ include as a special case the group $G$ itself. In this special
case, the ``(parabolic) character sheaves'' on $G$ are just the
usual character sheaves on $G$ introduced by Lusztig in \cite{L1}.
The varieties also include more or less as a special case the
boundary pieces of the De Concini-Procesi compactification of $G$
(where $G$ is adjoint).

\subsection*{0.2} We now review \cite{L3} and \cite{L4} in more detail.

For $Z_{J, y, D}$, there exists a finite partition into some smooth,
$G$-stable subvarieties which we call $G$-stable pieces. This
partition is based on some combinatorial result of B\'edard (see
2.2). The $G$-orbits on each piece are in one-to-one correspondence
with the ``twisted'' conjugacy classes of a certain (smaller)
reductive subgroup $L$. Furthermore, there is a natural equivalence
between the bounded derived category of $G$-equivariant,
constructible sheaves on that piece and the boundary derived
categories of $L$-equivariant (for the twisted conjugate action),
constructible sheaves on $L$.

To each character sheaf on $L$, one can associate a $G$-equivariant
simple perverse sheaf on the $G$-stable piece and call it a
character sheaf on the $G$-stable piece. This provides the ``local
picture'' of the theory of parabolic character sheaves. By imitating
the definition of character sheaves on the group, one can obtain
certain simple perverse sheaves on $Z_{J, y, D}$ and call them
character sheaves on $Z_{J, y, D}$. This provides the ``global
picture''.

Lusztig proved the following property:

Let $i$ be the inclusion of a $G$-stable piece to $Z_{J, y, D}$,
then

(1) for any character sheaf $C$ on $Z_{J, y, D}$, any perverse constituent
of $i^*(C)$ is a character sheaf on that piece;

(2) for any character sheaf $C$ on that piece, any perverse constituent of
$i_!(C)$ is a character sheaf on $Z_{J, y, D}$.

As a consequence, the character sheaves on $Z_{J, y, D}$ are just
the perverse extensions to $Z_{J, y, D}$ of the character sheaves on
$G$-stable pieces.

These results were also proved later in \cite{Sp2} and \cite{H2} in
some different way. In all these proofs, some inductive methods
based on B\'{e}dard's result were used. For more details, see the
introduction of \cite{H2}.

\subsection*{0.3} Let $\t$ be an involution on $G$ and $G^\t$ be the $\t$-fixed
point subgroup. In this paper, we study (under a mild assumption on
the characteristic of $\mathbf k$) a class of $G^\t$-equivariant
simple perverse sheaves on varieties $X_{J, \t}$ which we call
``(parabolic) character sheaves''. The varieties $X_{J, \t}$ are
defined in 3.5 and include the varieties $Z_{J, y, D}$ as some
special cases. They also include as a special case the strata of the
De Concini-Procesi compactification of the symmetric variety
$G/G^\t$ (hence the symmetric variety $G/G^\t$ itself). For more
details, see 3.30 and 3.31.

\subsection*{0.4} To achieve this goal, the first thing we need to
do is to find a partition of $X_{J, \t}$ which is analogous to the
partition of $Z_{J, y, D}$ into $G$-stable pieces. We call it the
partition of $X_{J, \t}$ into $G^\t$-stable pieces. However, there
is no results in our general setting that is analogous to
B\'{e}dard's result. Hence we need to find a different approach.

The idea is to relate the variety $X_{J, \t}$ to certain $Z_{J, y,
D}$. In the special case where $X_{J, \t}=G/G^\t$, we can identify
$G/G^\t$ with the identity component of $G^{\iota \circ \t}$, where
$\iota$ is the inverse map on $G$ (see \cite[3.3.0]{Gi}). This
result can be easily generalized. Namely, we can identify $X_{J,
\t}$ with certain irreducible component of $Z_{J, y, D}^{\iota \circ
\t}$, where $\iota$ is the ''inverse'' map on $Z_{J, y, D}$ (see
3.5).

Moreover, $\iota \circ \t$ maps an $G$-stable piece in $Z_{J, y, D}$
to another $G$-stable piece. This is what we will show in section 2.
Although this result is not needed to establish the partition of
$X_{J, \t}$, it serves as motivation for it. Namely, it suggests
that the $G^\t$-stable pieces might be the irreducible components of
the intersections of $X_{J, \t}$ with the $G$-stable pieces in
$Z_{J, y, D}$.

In fact, this is the right definition. (Certainly we need to show
that each connected component of the intersection is irreducible and
we need to know when the intersection is nonempty and what are the
components, etc.) Actually, in 3.21 we will define the $G^\t$-stable
pieces in an equivalent way that doesn't involve the $G$-stable
pieces in $Z_{J, y, D}$.

In section 3, we will also prove some results on the structure of
$G^\t$-stable pieces (see 3.26 and 3.28) and show that the
$G^\t$-orbits on each $G^\t$-stable piece are in one-to-one
correspondence with the $L^{\t_1}$-orbits on $L/L^{\t_2}$, where $L$
is a (smaller) reductive group and $\t_1, \t_2$ are two involutions
on $L$ (see 3.29).

\subsection*{0.5} Based on these results, we can establish a natural
equivalence between the bounded derived category of
$G^\t$-equivariant, constructible sheaves on that piece and the
boundary derived categories of $L^{\t_1}$-equivariant, constructible
sheaves on $L/L^{\t_2}$ (see 4.14).

Hence we obtain the ``local picture'' in the same way as in 0.2. The
``global picture'' is obtained by imitating Ginsburg's definition of
character sheaves on symmetric varieties in \cite{Gi}. There is also
a characterization of character sheaves using Ginsburg's
Harish-Chandra functor. This characterization will play an essential
role in our proofs of the properties analogous to those in 0.2.

Finally, in section 5, we generalize Lusztig's functors $e^{J'}_J$
and $f^J_{J'}$ and prove some properties.

\subsection*{Acknowledgement} We thank J. F. Thomsen
for carefully reading the first version of the manuscript and some
helpful suggestions. We also thank P. Deligne, G. Lusztig and D.
Vogan for some discussions.

\section{The $G$-stable pieces}

\subsection{} Let $B$ be a Borel
subgroup of $G$, $B^-$ be the opposite Borel subgroup and $T=B \cap
B^-$. Let $(\a_i)_{i \in I}$ be the set of simple roots. For $i \in
I$, we denote by $s_i$ the corresponding simple reflection. For any
element $w$ in the Weyl group $W=N(T)/T$, we use the same symbol $w$
for a representative of $w$ in $N(T)$. We denote by $\supp(w)
\subset I$ the set of simple roots whose associated simple
reflections occur in some (or equivalently, any) reduced
decomposition of $w$.

For $J \subset I$, let $W_J$ be the subgroup of $W$ generated by $J$
and $$W^J=\{w \in W; w=\min(w W_J)\}, \quad {}^J W=\{w \in W;
w=\min(W_J w)\}.$$ For $J, K \subset I$, we write $^J W^K$ for $^J W
\cap W^K$.

Let $\Phi^+$ be the set of positive roots in $G$ and $\Phi_J$ be the
set of roots in $L_J$. Set $\Phi_J^+=\Phi_J \cap \Phi^+$. Let $P_J
\supset B$ be the standard parabolic subgroup defined by $J$ and
$\cp_J$ be the set of parabolic subgroups conjugate to $P_J$. Let
$P^-_J \supset B^-$ be the opposite of $P_J$. Set $L_J=P_J \cap
P^-_J$. Then $L_J$ is a common Levi subgroup of $P_J$ and $P^-_J$.
Let $\pi_J: P_J \rightarrow L_J$ be the projection map.

For any parabolic subgroup $P$, we denote by $U_P$ its unipotent
radical. We simply write $U$ for $U_B$. For $J \subset I$, we simply
write $U_J$ for $U \cap L_J$ and $B_J$ for $B \cap L_J$.

For $J, K \subset I$, $P \in \cp_J, Q \in \cp_K$ and $u \in {}^J
W^K$, we write $\po(P, Q)=u$ if there exists $g \in G$, such that
$^g P=P_J$ and $^g Q={}^{ u} P_K$.

For any closed subgroup $H$ of $G$, we denote by $\Lie(H)$ the
corresponding Lie subalgebra and denote by $H_{\D}$ the image of the
diagonal embedding of $H$ in $G \times G$. For any subgroup $H$ and
$g \in G$, we write $^g H$ for $g H g\i$.

For each root $\a$, we denote by $u_{\a}$ the one-dimensional subspace in
$\Lie(G)$ that corresponds to $\a$.

If $\th$ be an automorphism on $G$ with $\th(T)=T$, then $\th$
induces a bijection on the set of roots and an automorphism on $W$.
By abuse of notation, we use the same symbol $\th$ for the induced
maps. We also use the same symbol $\th$ for the induced map on
$\Lie(G)$.

For a group,  we use $\iota$ for the inverse map. For an automorphism
$f$ on a variety $X$, we write $X^f$ for the fixed
point set.

\subsection{} Let $\hat G$ be an algebraic group with identity
component $G$ and $D$ be a fixed irreducible component of $\hat G$.
By the conjugation of Borel subgroups and maximal tori we may find
an element $g_D \in D$ such that $g_D$ normalizes $B$ and $T$. Set
$\d=\Ad(g_D)$.

Let $J, J' \subset I$ and $y \in {}^{J'} W^{\d(J)}$ be such that $y
\d(J)=J'$. For $P \in \cp_J$, $P' \in \cp_{J'}$, define $A_{y, D}(P,
P')=\{g \in D \mid \po(P', {}^g P)=y\}$. Define $$Z_{J, y, D}=\{(P,
P', \g); P \in \cp_J, P' \in \cp_{J'}, \g \in U_{P'} \backslash
A_{y, D}(P, P') /U_P\}$$ with $G \times G$-action defined by $(g_1,
g_2) (P, P', \g)=(^{g_2} P, {}^{g_1} P', g_1 \g g_2 \i)$.

By \cite[8.9]{L4}, $A_{y, D}(P, P')$ is a single $P', P$ double
coset. Thus $G \times G$ acts transitively on $Z_{J, y, D}$. Set
$$h_{J, y, D}=(P_J, {}^{\dot y \i} P_{J'}, U_{^{\dot y \i} P_{J'}}
g_D U_{P_J}) \in Z_{J, y, D}.$$ We call $h_{J, y, D}$ the base point
for the $G \times G$-action on $Z_{J, y, D}$. Now we may identify
$Z_{J, y, D}$ with $(G \times G)_{^{\dot y \i} P_{J'} \times P_J}
g_D L_J$ where $^{\dot y \i} P_{J'} \times P_J$ acts on $G \times G$
on the right and acts on $g_D L_J$ by $(p, p') \cdot g_D l=\bar p
g_D l \pi_J(p') \i$ for $p \in {}^{\dot y \i} P_{J'}$, $p' \in P_J$
and $l \in L_J$. Here $\bar p$ is the image of $p$ under the map
$^{\dot y \i} P_{J'} \rightarrow {}^{\dot y \i} P_{J'}/U_{^{\dot y
\i} P_{J'}} \cong L_{\d(J)}$.

\subsection{} For $w \in W^{\d(J)}$, set $I(J, w, \d)=\max\{K \subset J';
w \d(K)=K\}$ and $$Z_{J, y, D; w}=G_{\D} (B w, B) h_{J, y, D}.$$

The varieties $Z_{J, y, D; w}$ are called the $G$-stable pieces in
$Z_{J, y, D}$. They were introduced by Lusztig in \cite{L4}.

The following properties can be found in \cite[section 8]{L4} and \cite[section 1]{H1}.

(1) $Z_{J, y, D}=\bigsqcup_{w \in W^{\d(J)}} Z_{J, y, D; w}$.

(2) The map $G \times (P_{I(J, w, \d)}  w, P_{I(J, w, \d)}) h_{J, y,
D} \rightarrow Z_{J, y, D}$ defined by $(g, z) \mapsto (g, g) \cdot
z$ induces an isomorphism $$Z_{J, y, D; w} \cong G \times_{P_{I(J,
w, \d)}} (P_{I(J, w, \d)} w, P_{I(J, w, \d)}) h_{J, y, D},$$ where
the group $P_{I(J, w, \d)}$ acts on the right on $G$ and acts
diagonally on $(P_{I(J, w, \d)} w, P_{I(J, w, \d)}) h_{J, y, D}$.

(3) The map $$\bigl(U_{P_{I(J, w, \d)}} \cap {}^{w y \i} (U_J
U_{P^-_J}) \bigr) \times L_{I(J, w, \d)} w \rightarrow (P_{I(J, w,
\d)}  w, P_{I(J, w, \d)}) h_{J, y, D}$$ defined by $(u, z) \mapsto
(u z, 1) \cdot h_{J, y, D}$ is an isomorphism.

(4) The map $L_{I(J, w, \d)} w \rightarrow Z_{J, y, D}$ defined by
$z \mapsto (z, 1) \cdot h_{J, y, D}$ induces a bijection from the
set of $L_{I(J, w, \d)}$-orbits on $L_{I(J, w, \d)} w g_D$ to the
set of $G$-orbits on $Z_{J, y, D; w}$.

\section{The inverse map and the $G$-stable pieces}

\subsection{} Let $D$ be a connected component of $\hat G$. For the
connected component $D \i$ of $\hat G$, we choose $g_{D \i}=g_D \i$.
As in\cite[28.19]{L2}, we define the map $$\partial: Z_{J, y, D}
\rightarrow Z_{J', y \i, D \i}$$ by $\partial(P, Q, \g)=(Q, P, \g
\i)$. We call $\partial: Z_{J, y, D} \rightarrow Z_{J', y \i, D \i}$
the inverse map.

If $J=I$ and $y=1$, then $Z_{J, y, D}=D$ and $\partial$ is just the
restriction to $D$ of the inverse map on $\hat G$.

In Proposition 2.5, we will show that $\partial$ maps the $G$-stable
pieces in $Z_{J, 1, D}$ to the $G$-stable pieces in $Z_{\d(J), 1, D
\i}$.

\subsection{} In this subsection, we reformulate B\'edard's
description of $W^{\d(J)}$. The description below is slightly
different from \cite[2.2]{L3}. (In fact, the sequence $(J_n, w_n)_{n
\ge 0}$ below corresponds to the sequence $(J_n, w_n \i)_{n \ge 0}$
in {\it loc.cit.})

Let $\mathcal T(J, \d)$ be the set of all sequences $(J_n, w_n)_{n
\ge 0}$ where $J_n \subset I$ and $w_n \in W$ such that

(1) $J_0=J$,

(2) $J_n=J_{n-1} \cap \d \i (w_{n-1} \i J_{n-1})$ for $n \ge 1$,

(3) $w_n \in {}^{J_n} W^{\d(J_n)}$ and $w_n \in W_{J_{n-1}} w_{n-1}$
for $n \ge 1$.

For each sequence $(J_n, w_n) \in \mathcal T(J, \d)$, we have that
$J_m=J_{m+1}=\cdots$ and $w_m=w_{m+1}=\cdots$ for $m \gg 0$. By
\cite[2.4 \& 2.5]{L3}, $w_m \in W^{\d(J)}$ for all $m \ge 0$ and the
map $\mathcal T(J, \d) \rightarrow W^{\d(J)}$ defined by $(J_n,
w_n)_{n \ge 0} \mapsto w_m$ for $m \gg 0$ is a bijection. Moreover,
by \cite[1.4]{H1}, $J_m=I(J, w_m, \d)$ for $m \gg 0$.

\subsection{} To $(J_n, w_n)_{n \ge 0} \in \mathcal T(J, \d)$, we associate
a sequence $(K_n, v_n)_{n \ge 0}$ with $K_n \subset I$ and $v_n \in
W$. We set $$K_0=\d(J), \qquad v_0=w_0 \i.$$ Assume that $n \ge 1$
and that $K_m, v_m$ are already defined for $m<n$. Let
\begin{gather*} K_n=K_{n-1} \cap \d (v_{n-1} \i K_{n-1}), \\
v_n=\bigl(\d^n(v_0) \d^{n-1}(v_1) \cdots \d(v_{n-1}) \bigr) \i
\d^n(w_n \i) \bigl(\d^{n-1}(v_0) \d^{n-2}(v_1) \cdots v_{n-1}
\bigr).
\end{gather*}

This completes the inductive definition.

\begin{th1} We have that $(K_n, v_n)_{n \ge 0} \in \mathcal T(\d(J), \d \i)$.
\end{th1}

Proof. We show by induction on $n \ge 0$ that
\begin{equation*} K_n=\bigl(\d^n(v_0) \d^{n-1}(v_1) \cdots \d(v_{n-1}) \bigr) \i \d^{n+1}(J_n).\tag{a}
\end{equation*}

For $n=0$, (a) is clear. Assume now that $n>0$ and that the
statement holds when $n$ is replaced by $n-1$. Then
\begin{align*} K_{n-1} &=\bigl(\d^{n-1}(v_0) \d^{n-2}(v_1) \cdots \d(v_{n-2})
\bigr) \i \d^n(J_{n-1}) \\ &=\bigl(\d^n(v_0) \d^{n-1}(v_1) \cdots
\d(v_{n-1}) \bigr) \i \d^n(w_{n-1}) \i \d^{n+1}(J_n) \end{align*}
and
\begin{align*} K_n &=K_{n-1} \cap \d (v_{n-1} \i K_{n-1}) \\ &=\bigl(\d^n(v_0) \d^{n-1}(v_1) \cdots \d(v_{n-1}) \bigr) \i
\d^n(w_{n-1} \i J_{n-1}) \\ & \qquad \cap \bigl(\d^n(v_0) \d^{n-1}(v_1) \cdots \d(v_{n-1}) \bigr) \i \d^{n+1}(J_{n-1}) \\
&=\bigl(\d^n(v_0) \d^{n-1}(v_1) \cdots \d(v_{n-1}) \bigr) \i \d^{n+1}(\d \i (w_{n-1} \i J_{n-1}) \cap J_{n-1}) \\
&=\bigl(\d^n(v_0) \d^{n-1}(v_1) \cdots \d(v_{n-1}) \bigr) \i \d^{n+1}(J_n).
\end{align*}

(a) is proved.

By definition, $v_0 \in {}^{K_0} W^{\d \i(K_0)}$.

Let $n \ge 1$. By definition, $w_n \in {}^{J_n} W^{\d(J_n)}$,
$w_{n-1} \in {}^{J_{n-1}} W^{\d(J_{n-1})}$ and $w_n w_{n-1} \i \in
W_{J_{n-1}}$. Thus \begin{align*} \tag{b} w_n \d(\Phi_{J_n}^+) &=w_n
\d(\Phi_{J_n}) \cap \Phi^+ \subset w_n w_{n-1} \i(\Phi_{J_{n-1}})
\cap \Phi^+ \\ &=\Phi_{J_{n-1}} \cap \Phi^+=\Phi_{J_{n-1}}^+, \\
\tag{c} w_n \i (\Phi_{J_n}^+) &=w_n \i(\Phi_{J_n}) \cap \Phi^+
\subset w_n \i(\Phi_{J_{n-1}}) \cap \Phi^+ \\ &=w_{n-1}
\i(\Phi_{J_{n-1}}) \cap \Phi^+=w_{n-1} \i(\Phi_{J_{n-1}}^+).
\end{align*}

From (a) and the definition of $v_n$, we deduce that
\begin{gather*} \tag{d} v_n \d \i(\Phi_{K_n}^+)=\bigl(\d^n(v_0) \d^{n-1}(v_1)
\cdots \d(v_{n-1}) \bigr) \i \d^n(w_n) \i \d^n(\Phi_{J_n}^+), \\
\tag{e} v_n \i(\Phi_{K_n}^+)=\bigl(\d^{n-1}(v_0) \d^{n-2}(v_1)
\cdots v_{n-1} \bigr) \i \d^n(w_n) \d^{n+1}(\Phi_{J_n}^+).
\end{gather*}

From (c) and (d), we see that \begin{align*} v_n \d \i(\Phi_{K_n}^+)
& \subset \bigl(\d^n(v_0) \d^{n-1}(v_1) \cdots \d(v_{n-1})
\bigr) \i \d^n(w_{n-1}) \i \d^n(\Phi_{J_{n-1}}^+) \\
&=\bigl(\d^{n-1}(v_0) \d^{n-2}(v_1) \cdots \d(v_{n-2}) \bigr) \i
\d^n(\Phi_{J_{n-1}}^+)=\Phi_{K_{n-1}}^+. \end{align*}

From (b) and (e), we see that $$v_n \i(\Phi_{K_n}^+) \subset
\bigl(\d^{n-1}(v_0) \d^{n-2}(v_1) \cdots v_{n-1} \bigr) \i
\d^n(\Phi_{J_{n-1}}^+)=v_{n-1} \i \Phi_{K_{n-1}}^+.$$

Using induction method, we deduce $v_n \i(\Phi_{K_n}^+)
\subset \Phi^+$ for all $n \ge 0$.

Therefore $v_n \in {}^{K_n} W^{\d\i(K_n)}$.

Let $n \ge 1$. By definition, \begin{align*} \d(v_{n-1})= &
\bigl(\d^n(v_0) \d^{n-1}(v_1) \cdots \d^2(v_{n-2}) \bigr) \i
\d^n(w_{n-1}) \i \\ & \bigl(\d^{n-1}(v_0) \d^{n-2}(v_1) \cdots
\d(v_{n-2})\bigr).\end{align*}

Thus
\begin{align*} \tag{f} v_n=& \bigl(\d^{n-1}(v_0) \d^{n-2}(v_1) \cdots \d(v_{n-2})
\bigr) \i \d^n(w_{n-1} w_n \i) \\ & \bigl(\d^{n-1}(v_0)
\d^{n-2}(v_1) \cdots v_{n-1} \bigr).\end{align*}

Notice that $w_{n-1} w_n \i \in W_{J_{n-1}}$. By (a), $v_n \in W_{K_{n-1}} v_{n-1}$. \qed

\begin{th2} Define a map $\e_{J, \d}: W^{\d(J)} \rightarrow W^J$
by sending the element $w \in W^{\d(J)}$ corresponding to $(J_n,
w_n)_{n \ge 0}$ to the element in $W^J$ corresponding to $(K_n,
v_n)_{n \ge 0}$. Then $\partial(Z_{J, 1, D; w})=Z_{\d(J), 1, D\i;
\e_{J, \d}(w)}$.
\end{th2}

Proof. For $n \ge 0$, set \begin{align*} x_n= &\bigl(\d^n(v_0)
\d^{n-1}(v_1) \cdots \d(v_{n-1}) \bigr) \i \d^{n+1} (w_n w \i) \\ &
\bigl(\d^n(v_0) \d^{n-1}(v_1) \cdots \d(v_{n-1}) \bigr).\end{align*}
Since $w_n w \i \in W_{J_n}$, then $x_n \in W_{K_n}$. Therefore
$l(x_n v_n)=l(x_n)+l(v_n)$ and $l(v_n \d \i(x_n))=l(v_n)+l(x_n)$. By
2.4(f),
\begin{align*} & x_n v_n=\bigl(\d^n(v_0) \d^{n-1}(v_1) \cdots \d(v_{n-1})
\bigr) \i \d^{n+1} (w_n w \i) \bigl(\d^n(v_0) \d^{n-1}(v_1) \cdots
v_n \bigr) \\ & =v_{n+1} \bigl(\d^n(v_0) \d^{n-1}(v_1) \cdots v_n
\bigr) \i \d^{n+1} (w_{n+1} w \i) \bigl(\d^n(v_0) \d^{n-1}(v_1)
\cdots v_n \bigr) \\& =v_{n+1} \d \i(x_{n+1}).
\end{align*}

By definition, \begin{align*} \partial(Z_{J, 1, D; w}) &=\partial(G_{\D} (B w, 1)
\cdot h_{J, 1, D})=G_{\D} (1, B  w) \cdot h_{\d(J), 1, D \i} \\
&=G_{\D} ( w \i B, 1) \cdot h_{\d(J), 1, D \i}=G_{\D} B_{\D} ( w \i
B, 1) \cdot h_{\d(J), 1, D \i}
\\ &=G_{\D} (B  v_0 \d \i(x_0) B, 1) \cdot h_{\d(J), 1, D \i}. \end{align*}

For $n \ge 0$, \begin{align*} & G_{\D} (B  v_n \d \i( x_n) B, 1)
\cdot h_{\d(J), 1, D \i}=G_{\D} (B  v_n B \d \i( x_n) B, 1) \cdot
h_{\d(J), 1, D \i}
\\ &=G_{\D} (B  v_n B, B  x_n \i B) \cdot h_{\d(J), 1, D \i}=G_{\D}
(B  x_n B  v_n B, 1) \cdot h_{\d(J), 1, D \i} \\ &=G_{\D} (B
v_{n+1} \d \i( x_{n+1}) B, 1) \cdot h_{\d(J), 1, D \i}. \end{align*}

Therefore $\partial(Z_{J, 1, D; w})=G_{\D} (B v_n \d \i(x_n), 1)
\cdot h_{\d(J), 1, D\i}$ for all $n \ge 0$. In particular,
$$\partial(Z_{J, 1, D; w})=G_{\D} (B  \e_{J, \d}(w), B) \cdot
h_{\d(J), 1, D\i}=Z_{\d(J), 1, D \i; \e_{J, \d}(w)}.$$ \qed

\

Notice that $\partial \circ \partial (Z_{J, 1, D; w})=Z_{J, 1, D; w}$.
We have the following consequence.

\begin{th3} The map $\e_{\d(J), \d \i} \circ \e_{J, \d}: W^{\d(J)}
\rightarrow W^{\d(J)}$ is the identity map.
\end{th3}

\

The following corollary gives another characterization of the map $\e_{J, \d}$.

\begin{th4} For each $w \in W^{\d(J)}$, there exists a unique
element in $W^J$ which is of the form $\d(x) \i w \i x$ for some $x
\in W_J$. This element is just $\e_{J, \d}(w)$.
\end{th4}

Proof. The existence of the element follows from the proof of
Proposition 2.5. We prove the uniqueness. Assume that $\d(x) \i w \i
x \in W^J$ for some $x \in W_J$. Then $$({\d(x)} \i w \i  x, 1)
\cdot h_{\d(J), 1, D\i} \in Z_{\d(J), 1, D \i; \d(x) \i w \i x}.$$
On the other hand, \begin{align*} ({\d(x)} \i w \i  x, 1) \cdot
h_{\d(J), 1, D\i} &=( {\d(x)} \i, {\d(x)} \i) ( w \i  x, {\d(x)})
\cdot h_{\d(J), 1, D\i} \\ & \in G_{\D} ( w \i T, 1) \cdot h_{\d(J),
1, D\i} \subset \partial(Z_{J, 1, D; w})\\ &=Z_{\d(J), 1, D \i;
\e_{J, \d}(w)}. \end{align*}

Hence $\d(x) \i w \i x=\e_{J, \d}(w)$. \qed

\

In general, $\partial$ doesn't map a $G$-stable piece in $Z_{J, y,
D}$ to a $G$-stable piece in $Z_{J', y \i, D \i}$. However, we have
a modified version which will be stated in the end of this section.

\begin{th5} Let $x \in W$ and $L$ be a Levi of some parabolic subgroup
of $L_J$ with $^{x g_D} L=L$. Then there exists $w \in W^{\d(J)}$,
such that $$(L x, 1) \cdot h_{J, y, D} \subset Z_{J, y, D; w}.$$
\end{th5}

Proof. We define by induction on $n$ a sequence $(J_n, w_n, u_n)_{n
\ge 0}$ as follows.

We write $x$ as $a \d(b)$ for $a \in W^{\d(J)}$ and $b \in W_J$ and set
$$J_0=J, w_0=\min(W_J a), u_0=b a w_0 \i.$$ Assume that $n>0$ and that
$J_{n-1}, w_{n-1}, u_{n-1}$ are defined. Let $$J_n=J_{n-1} \cap \d
\i(w_{n-1} \i J_{n-1}).$$ Write $u_{n-1} w_{n-1}$ as $w'_n \d(u'_n)$
for $w'_n \in W^{\d(J_n)}$ and $u'_n \in W_{J_n}$. Set
$$w_n=\min(W_{J_n} w'_n), u_n=u'_n w'_n w_n \i.$$ This completes the
inductive definition.

We show that \begin{equation*} w_n \in {}^{J_n} W^{\d(J_n)} \quad \text{ for } n \ge 0. \tag{a}\end{equation*}

Let $y=\min(W_{J_n} w_n W_{\d(J_n)})$. Now $w_n' \in W_{J_n} y
W_{\d(J_n)}$ and $w_n' \in W^{\d(J_n)}$. By \cite[2.1(b)]{L3}, $w'_n
\in W_{J_n} y$. Notice that $w'_n \in W_{J_n} w_n$ and $w_n, y \in
{}^{J_n} W$. Then $w_n=y$ and (a) follows.
\begin{equation*} w_n \in W_{J_{n-1}} w_{n-1} \quad \text{ for } n \ge 1. \tag{b}\end{equation*}

By definition, $w_n w_{n-1} \i \in W_{J_n} w_n' w_{n-1} \i \subset
W_{J_n} u_{n-1} w_{n-1} W_{\d(J_n)} w_{n-1} \i$. Since $u_{n-1} \in
W_{J_{n-1}}$ and $J_n, w_{n-1} \d(J_n) \subset J_{n-1}$, we have
$w_n w_{n-1} \i \in W_{J_{n-1}}$ and (b) follows.
\begin{equation*} (J_n, w_n)_{n \ge 0} \in \ct(J,
\d).\tag{c}\end{equation*} This follows from (a) and (b).

For $n \ge 0$, set $L_n={}^{(u'_n u'_{n-1} \cdots u'_1) b} L$. We
show by induction on $n \ge 0$ that \begin{gather*} {}^{u_n w_n g_D}
L_n=L_n, \tag{d} \\ (L x, 1) \cdot h_{J, y, D} \subset G_{\D} (L_n
u_n w_n, 1) \cdot h_{J, y, D}. \tag{e}\end{gather*}

For $n=0$, ${}^{u_0 w_0 g_D} L_0={}^{b a \d(b) g_D} L={}^b ({}^{x
g_D} L)=L_0$ and \begin{align*} (L x, 1) \cdot h_{J, y, D} &=(L a, b
\i) \cdot
h_{J, y, D} \subset G_{\D} (b L a, 1) \cdot h_{J, y, D}\\
&=G_{\D}(L_0 u_0 w_0, 1) \cdot h_{J, y, D}. \end{align*}

Assume now that $n>0$ and that (d) and (e) hold when $n$ is replaced
by $n-1$. Then \begin{align*} {}^{u_n w_n g_D} L_n &={}^{u'_n w'_n
\d(u_n') g_D} L_{n-1}={}^{u'_n} ({}^{u_{n-1} w_{n-1} g_D}
L_{n-1})={}^{u'_n} L_{n-1}=L_n, \\ (L x, 1) \cdot h_{J, y, D} &
\subset G_{\D} (L_{n-1} w'_n, (u'_n) \i) \cdot h_{J, y, D}=G_{\D}
(u'_n L_{n-1} w'_n, 1) \cdot h_{J, y, D} \\ &=G_{\D} (L_n u_n w_n,
1) \cdot h_{J, y, D}.
\end{align*}
Thus (d) and (e) are proved.

By (c), there exists $m>0$ such that $w_m \in W^{\d(J)}$, $u_m \in
W_{J_m}$ and $w_m \d(J_m)=J_m$. By \cite[Lemma 7.3]{St}, $u_m=l_1
l_2 (w_m g_D) l_1 \i (w_m g_D) \i$ for $l_1 \in L_{J_m}$ and $l_2
\in L_{J_m} \cap B$. Set $L'={}^{l_1 \i} L_m$. Then $^{l_2 w_m g_D}
L'=L'$ and
\begin{align*} & (L x, 1) \cdot h_{J, y, D} \subset G_{\D} (L_m u_m
w_m, 1) \cdot h_{J, y, D}=G_{\D} (L_m l_1 l_2 w_m, l_1) \cdot h_{J,
y, D} \\ &=G_{\D} (L' l_2 w_m, 1) \cdot h_{J, y, D}=G_{\D}
\bigl(L'_{\D} ((L' \cap B) l_2 w_m, 1) \cdot h_{J, y, D}\bigr) \\ &
\subset G_{\D} (B w_m, 1) \cdot h_{J, y, D}=Z_{J, y, D; w_m}.
\end{align*}

\begin{th6} Let $\s$ be an involution on $\hat G$
with $\s(g_D)=g_D \i$ and $\s(P_J)={}^{y \i} P_{J'}$. Define the map
$\s: Z_{J, y, D} \rightarrow Z_{J', y \i, D \i}$ by $$\s(P, Q,
\g)=(\s(P), \s(Q), \s(\g)).$$ Then there exists a map $\wp:
W^{\d(J)} \rightarrow W^{\d(J)}$ such that $$\partial \circ \s
(Z_{J, y, D; w})=Z_{J, y, D; \wp(w)}$$ for all $w \in W^{\d(J)}$. In
particular, $\wp \circ \wp=id$.
\end{th6}

Proof. We have that $\partial \circ \s (Z_{J, y, D})=Z_{J, y, D}$
and $\partial \circ \s(h_{J, y, D})=h_{J, y, D}$. Then
\begin{align*}
\partial \circ \s (Z_{J, y, D; w}) &=\partial \circ \s \bigl(G_{\D}
(w L_{\d(I(J, w, \d))}), 1) \cdot h_{J, y, \d} \bigr) \\ &=G_{\D}
\bigl(\s(L_{\d(I(J, w, \d))}) \s(w) \i, 1 \bigr) \cdot h_{J, y,
\d}.\end{align*} Notice that $$^{\s(w) \i g_D} \s(L_{\d(I(J, w,
\d))})=\s({}^{w \i g_D \i} L_{\d(I(J, w, \d))})=\s(L_{\d(I(J, w,
\d))}).$$ By Lemma 2.8, there exists a map $\wp: W^{\d(J)}
\rightarrow W^{\d(J)}$ such that such that $\partial \circ \s (Z_{J,
y, D; w}) \subset Z_{J, y, D; \wp(w)}$.

Then $\partial \circ \s (Z_{J, y, D; \wp(w)}) \subset Z_{J, y, D;
\wp \circ \wp(w)}$. Since $(\partial \circ \s) \circ (\partial \circ
\s)=id$, $\partial \circ \s (Z_{J, y, D; \wp(w)}) \subset Z_{J, y,
D; w}$. Therefore $\partial \circ \s (Z_{J, y, D; w})=Z_{J, y, D;
\wp(w)}$ and $\wp \circ \wp=id$. \qed

\section{The $G^\s$-stable pieces}

\subsection{} From now on, we assume that the characteristic of
$\mathbf k$ is $0$ or sufficiently large. Let $\s$ and $\t$ be involutions on $G$ with $\s(T)=\t(T)=T$
and $\s \t(B)=B$. Set $\hat G=G \ltimes <\s \t>$. Then $\s$, $\t$
acts on $\hat G$ by \begin{align*} \s(g, (\s \t)^n) &=(\s(g), (\t
\s)^n)=(\s(g), (\s \t)^{-n}), \\ \t(g, (\s \t)^n) &=(\t(g), (\s
\t)^{-n}). \end{align*}

Let $D=(G, \s \t)$. Then $D$ is a connected component of $\hat G$
and $\s(D)=D \i$. Set $g_D=(1, \s \t)$. Then $\s(g_D)=g_D \i$. As in
1.2, we use the symbol $\d$ for the induced maps of $\s \t$ on
$\Phi$, $I$ and $W$.

\subsection{} Let $J \subset I$ with $\t(L_J)=L_J$. Then $\s(P_J),
P_{\d(J)}=\s \t(P_J)$ have a common Levi. Let $J' \subset I$ and $y
\in {}^{J'} W^{\d(J)}$ be such that $\s(P_J)={}^{y \i} P_{J'}$. By
\cite[2.8.7]{C}, we have that $y \d(J)=J'$.

Recall that $L_J^{\iota \circ \t}=\{l \in L_J; \t(l)=l \i\}$. Define
the $P_J$-action on $G \times L_J^{\iota \circ \t}$ by $p \cdot (g,
l)=(g p \i, \pi_J(p) l \t(\pi_J(p)) \i)$. Let $G \times_{P_J}
L_J^{\iota \circ \t}$ be the quotient space. In Proposition 3.4, we
will show that $G \times_{P_J} L_J^{\iota \circ \t}$ can be
identified with certain subvariety of $Z_{J, y, D}$. Before doing
that, let us recall the following result (see \cite[page 26, lemma
4]{SL}).

\begin{th7} Let $H$ be a closed subgroup of $G$ and $\Phi:
X \rightarrow G/H$ be a $G$-equivariant morphism from the
$G$-variety $X$ to the homogeneous space $G/H$. Let $E \subset X$ be
the fiber $\Phi \i(H)$. Then $E$ is stabilized by $H$ and the map
$\Psi: G \times_H E \rightarrow X$ sending $(g, e)$ to $g \cdot e$
defines an isomorphism of $G$-varieties.
\end{th7}

\

Now we prove the following result.

\begin{th8} The map $G \times L_J^{\iota \circ \t}
\rightarrow Z_{J, y, D}$ defined by $(g, l) \mapsto (\s(g), g l)
\cdot h_{J, y, D}$ induces an isomorphism $$G \times_{P_J}
L_J^{\iota \circ \t} \cong Z_{J, y, D}^{\iota \circ \s}.$$ In
particular, $Z_{J, y, D}^{\iota \circ \s}$ is irreducible if and
only if $L_J^{\iota \circ \t}$ is irreducible.
\end{th8}

\begin{rmk} This result is inspired by \cite[2.3]{Sp1}.
\end{rmk}

Proof. Define a $G$-action on $Z_{J, y, D}$ by $g \cdot z=(\s(g), g)
\cdot z$.  Then $Z_{J, D}^{\iota \circ \s}$ is stable under the
action. The $G$-equivariant morphism $Z_{J, y, D} \rightarrow G/P_J$
defined by $(P, Q, \r) \mapsto P$ induces a $G$-equivariant morphism
$$\Phi: Z_{J, y, D}^{\iota \circ \s} \rightarrow G/P_J.$$

If $(g', g) \cdot h_{J, y, D} \in \Phi \i(P_J)$, then $g=l u$ for
some $l \in L_J$ and $u \in U_{P_J}$ and $(g', g) \cdot h_{J, y,
D}=(\s(g), \s(g')) \cdot h_{J, y, D}=(1, \s(g') \t(l) \i) \cdot
h_{J, y, D}$. Thus $\s(g') \in P_J$. Write $\s(g')$ as $l' u'$ for
$l' \in L_J$ and $u' \in U_{P_J}$. Then
\begin{align*} & (g', g) \cdot h_{J, y, D}=(\s(l'), l)
\cdot h_{J, y, D}=(1, l \t(l') \i) \cdot h_{J, y, D}, \\ & (\s(g),
\s(g')) \cdot h_{J, y, D}=(\s(l), l') \cdot h_{J, y, D}=(1, l' \t(l)
\i) \cdot h_{J, y, D}.\end{align*}

So $l \t(l') \i=\t(l' \t(l) \i) \i$ and $\Phi \i(P_J)=(1, L_J^{\iota
\circ \t}) \cdot h_{J, y, D}$. Now the proposition follows from
Lemma 3.3. \qed

\subsection*{3.5} By \cite[Lemma 3.3.0]{Gi}, the ``Lang map'' $l
\mapsto l \t(l) \i$ gives rise to a $L_J$-equivariant isomorphism
from $L_J/L_J^{\t}$ to the identity component of $L_J^{\iota \circ
\t}$. Now define the $P_J$ action on $G \times L_J/L_J^{\t}$ by $p
\cdot (g, z)=(g p \i, \pi_J(p) z)$. Let $$X_{J, \t}=G \times_{P_J}
L_J/L_J^{\t}$$ be the quotient space. Then by the identification of
$G \times_{P_J} L_J^{\iota \circ \t}$ with $Z_{J, y, D}^{\iota \circ
\s}$, we may identify $X_{J, \t}$ with the irreducible component of
$Z_{J, y, D}^{\iota \circ \s}$ that contains $h_{J, y, D}$.

We can also naturally identify $X_{J, \t}$ with $G/U_{P_J} L_J^\t$.
In particular, $G$ acts transitively on $X_{J, \t}$.

\subsection*{3.6. Examples} (1). We write $\mathbf G=G \times G$,
$\mathbf T=T \times T$, $\mathbf B=B \times B$. Denote by $\s$ the
permutation involution $(g, g') \rightarrow (g', g)$ of $\mathbf G$.
Let $J \subset I$ and $\mathbf J=(J, J)$. Then $\mathbf X_{\mathbf
J, \s}=Z_{J, 1, G}$ and the $\mathbf G^{\s}$-action on $\mathbf
X_{\mathbf J, \s}$ is just the $G_{\D}$-action on $Z_{J, 1, G}$.

(2). Let $\s$ be an involution on $G$. Then $X_{I, \s}=G/G^{\s}$ is a symmetric space.

\begin{th9} Let $\th$ be an involution on $G$ such that $\th(T)=T$.
Then for each $G^{\th} \times B$-orbit $v$ on $G$, there exists $x
\in v$ such that $x \i \th(x) \in N(T)$.
\end{th9}

See \cite[Theorem 1.3]{RS} for the special case where $\th(B)=B$ and
\cite[3.5]{Sp1} for the general case.

\subsection*{3.8} Let $\th$ be an involution on $G$ with $\th(T)=T$. Let $J \subset I$
with $\th(L_J)=L_J$. Define $$\cj_{J, \th}=\{w \in W; g \i \th(g)
\in B_J w \, \th(B_J) \text{ for some } g \in L_J\}.$$ Set $$\cw(J,
\s, \t)=\{w \in W^{\d(J)}; w \d(u) \in \cj_{I, \s} \text{ for some }
u \in \cj_{J, \t}\}.$$

\begin{th10} Let $w \in W^{\d(J)}$. Then $Z_{J, y, D; w}
\cap X_{J, \t}\neq \varnothing$ if and only if $w \in \cw(J, \s, \t)$.
\end{th10}

\begin{rmk} If $Z_{J, y, D; w} \cap X_{J, \t}\neq \varnothing$, then
$Z_{J, y, D; w} \cap Z_{J, y, D}^{\iota \circ \s} \neq \varnothing$.
In other words, $\iota \circ \s(Z_{J, y, D; w}) \cap Z_{J, y, D; w}
\neq \varnothing$. By Proposition 2.9, $\iota \circ \s(Z_{J, y, D;
w})=Z_{J, y, D; w}$.
\end{rmk}

Proof. Let $g \in G$ and $b \in B$ with $(g b w, g) \cdot h_{J, D} \in X_{J, \t}$. Then
$$h_{J, y, D}=(w \i b \i g \i \s(g), g \i \s(g) \s(b) \s(w)) \cdot h_{J, y, D}.$$
So $w \i b \i g \i \s(g) \in \s(P_J)$ and $g \i \s(g) \in B w
\s(P_J)=B w L_{\d(J)} \s(U_{P_J})$. In other words, $g \i \s(g) \in
B w \d(u) \s(B)$ for some $u \in W_J$. By definition, $w \d(u) \in
\cj_{I, \s}$.

By Proposition 3.7, we may write $g$ as $k x b'$ for $k \in G^{\s}$,
$x \in G$ with $x \i \s(x) \in w \d(u) T$ and $b' \in B$. Notice
that \begin{align*}(k x B w, k x B) \cdot h_{J, y, D} &=(k x w
({}^{w \i} B), k x B) \cdot h_{J, y, D} \\ &=(k \s(x) \d(u) \i
({}^{w \i} B), k x B) \cdot h_{J, y, D}. \end{align*} Thus $ \d(u)
\i ({}^{w \i} B), B) \cdot h_{J, y, D} \cap X_{J, \t} \neq
\varnothing$.

In other words, there exists $b_1 \in {}^{w \i} B$, $b_2 \in B$ and
$p \in G$, such that $(\d(u) \i b_1, b_2) \cdot h_{J, y, D}=(\s(p),
p) \cdot h_{J, y, D}$. Hence $b_2 \i p \in P_J$ and $p \in P_J$.
Therefore $b_1 \in \d(u) \s(p) \s(P_J)=\s(P_J)$. Write $b_1$ as
$b_1=l b_1'$ for $l \in L_{\d(J)}$ and $b_1' \in U_{\s(P_J)}$. Then
it is easy to see that $l \in {}^{w \i} B$. Since $w \in W^{\d(J)}$,
$l \in L_{\d(J)} \cap B$. Therefore \begin{align*} & (\d(u) \i b_1,
b_2) \cdot h_{J, y, D}=(1, b_2 g_D \i l \i \d(u) g_D) \cdot h_{J, y,
D} \in (1, B u) \cdot h_{J, y, D} \\ & (\s(p), p) \cdot h_{J, y,
D}=(1, \pi_J(p) \t(\pi_J(p)) \i) \cdot h_{J, D}.\end{align*}

Therefore, $\pi_J(p) \t(\pi_J(p)) \i \in B u \t(B)$ and $u \in
\cj_{J, \t}$. Hence $w \in \cw(J, \s, \t)$.

On the other hand, assume that $u \in W_J$ with $w \d(u) \in \cj_{I,
\s}$ and $u \in \cj_{J, \t}$. Then by Proposition 3.7, there exists
$x \in G$ with $x \i \s(x) \in w {\d(u)} T$ and $y \in L_J$ with $y
\i \t(y) \in u T$. Then \begin{align*} (\s(x y \i), x y \i) \cdot
h_{J, y, D} &=(x, x) (x \i \s(x), y \i \t(y)) \cdot h_{J, y, D} \\ &
\in (x, x) (w \d(u) T, u T) \cdot h_{J, y, D} \\ &=(x, x) (w T, 1)
\cdot h_{J, y, D} \subset Z_{J, y, D; w}. \end{align*}

Thus $Z_{J, y, D; w} \cap X_{J, \t} \neq \varnothing$. \qed

\

We can see that for $w \in \cw(J, \s, \t)$, the elements $u \in W_J$
such that $w \d(u) \in \ci_{I, \s}$ and $u \in \ci_{J, \t}$ may not
be unique. We will study $\cw(J, \s, \t)$ in more detail and show
that there exists a ``distinguished'' element $u$ for each $w \in
\cw(J, \s, \t)$.

\begin{th11} Let $x, w \in W^{\d(J)}$, $u \in W_J$ and $v
\in W_{\d(J)}$. If $\supp(v)=\d(\supp(u))$ and $w v=u x$, then $u
\in W_{I(J, w, \d)}$.
\end{th11}

Proof. Let $(J_n, w_n)_{n \ge 0}$ be the element in $\ct(J, \d)$
that corresponds to $w$. By 2.2, it suffices to prove that $u \in
W_{J_n}$ for $n \ge 0$.

We argue by induction on $n$. For $n=0$ this is clear. Assume now
that $n>0$ and that (a) holds when $n$ is replaced by $n-1$. Write
$v$ as $a b$ for $a \in W_{\d(J_n)}$ and $b \in {}^{\d(J_n)} W \cap
W_{\d(J_{n-1})}$ and $x$ as $x=c d$ for $c \in W_{J_{n-1}}$ and $d
\in {}^{J_{n-1}} W^{\d(J)}$. Then $w_{n-1} v=(w_{n-1} a w_{n-1} \i)
w_{n-1} b$, where $w_{n-1} a w_{n-1} \i \in W_{J_{n-1}}$ and
$w_{n-1} b \in {}^{J_{n-1}} W$. Thus $w v \in W_{J_{n-1}} w_{n-1}
v=W_{J_{n-1}} w_{n-1} b$ and $u x \in W_{J_{n-1}} d$. Then $w_{n-1}
b=d$. Since $d, w_{n-1} \in W^{\d(J)}$, we see that $b=1$ and $v \in
W_{\d(J_n)}$. By our assumption, $u \in W_{J_n}$. \qed

\begin{th12} Let $w \in W^{\d(J)}$ and $u \in W_J$. If
$\s(w)=\d(u) \i w \i u$, then $u \t(\Phi_{I(J, w, \d)})=\Phi_{I(J,
w, \d)}$.
\end{th12}

Proof. Since $w \d(I(J, w, \d))=I(J, w, \d)$, there exists a
bijection
$$\rho: I(J, w, \d) \rightarrow I(J, w, \d)$$ such that $w
\a_{\d(j)}=\a_{\rho(j)}$ for each $k \in I(J, w, \d)$.

Applying $\s$ on both sides, we have that $\d(u) \i w \i u \s
(\a_{\d(j)})=\s (\a_{\rho(j)})$. Hence $$u \t( \a_j)=u \s
(\a_{\d(j)})=w \d(u) \s (\a_{\rho(j)})=w \d \bigl(u \t(
\a_{\rho(j)}) \bigr).$$ Then
\begin{align*} \sum_{j \in I(J, w, \d), u \t( \a_j)>0} u \t( \a_j)
&=w \d(\sum_{j \in I(J, w, \d), u \t( \a_j)>0} u \t( \a_j)), \\
\sum_{j \in I(J, w, \d), u \t( \a_j)<0} u \t( \a_j) &=w \d(\sum_{j
\in I(J, w, \d), u \t( \a_j)<0} u \t( \a_j)).
\end{align*}
Assume that \begin{align*} \sum_{j \in I(J, w, \d), u \t( \a_j)>0} u
\t( \a_j) &=\sum_{i \in J} a_i \a_i, \\ \sum_{j \in I(J, w, \d), u
\t( \a_j)<0} u \t(\a_j) &=-\sum_{i \in J} b_i \a_i. \end{align*}
Since $w \in W^{\d(J)}$, then $w \d\{i \in J; a_i \neq 0\}=\{i \in
J; a_i \neq 0\}$ and $w \d \{i \in J; b_i \neq 0\}=\{i \in J; b_i
\neq 0\}$.

Moreover, for each $j \in I(J, w, \d)$, $u \t( \a_j)$ is a linear
combination of $\a_i$ with $a_i \neq 0$ or $b_i \neq 0$. Thus the
vector space spanned by $u \t( \a_j)$ for $j \in I(J, w, \d)$ is a
subspace of the vector space spanned by $\a_i$ with  $a_i \neq 0$ or
$b_i \neq 0$ and the cardinality of $\{i \in J; a_i \neq 0, \text{
or } b_i \neq 0\}$ is larger than or equal to the cardinality of
$I(J, w, \d)$.

By the definition of $I(J, w, \d)$, we have that $I(J, w, \d)=\{i
\in J; a_i \neq 0, \text{ or } b_i \neq 0\}$. Therefore, the vector
space spanned by $u \t( \a_j)$ for $j \in I(J, w, \d)$ equals the
vector space spanned by $\a_k$ for $k \in I(J, w, \d)$. Hence $u
\t(\Phi_{I(J, w, \d)})=\Phi_{I(J, w, \d)}$. The lemma is proved.
\qed

\begin{th13} Let $\th$ be an involution on $G$ with
$\th(T)=T$. Let $K \subset I$, $w \in W$ with $w
\th(\Phi_K^+)=\Phi_K^+$. If $g \i \th(g) \in P_K w \th(P_K)$. Then
there exists $l \in L_K$, such that $(g l) \i \th(g l) \in B w
\th(B)$.
\end{th13}

\begin{rmk} This result is due to J. F. Thomsen by private communication.
\end{rmk}

Proof. Set $\th'=\Ad(w) \circ \th$. Then $\th'$ is an automorphism
on $G$. By our assumption, $\th'(L_K)=L_K$ and $\th'(B_K)=B_K$. We
have that $g \i \th(g) \in U_{P_K} l_1 w \th(U_{P_K})$ for some $l_1
\in L_K$.

By \cite[Lemma 7.3]{St}, the map $L_K \times B_K \rightarrow L_K$
defined by $(l, b) \mapsto l \i b \th'(l)$ is surjective. Thus $l \i
l_1 \th'(l) \in B_K$ for some $l \in L_K$. Hence $(g l) \i \th(g l)
\in U_{P_K} B_K w \th(U_{P_K}) \subset B w \th(B)$. The lemma is
proved. \qed

\begin{th14} For $w \in \cw(J, \s, \t)$, there exists a
unique element $u \in W_J$ such that $u \t( \Phi_{I(J, w,
\d)}^+)=\Phi_{I(J, w, \d)}^+$, $w \d(u) \in \cj_{I, \s}$ and $u \in
\cj_{J, \t}$.
\end{th14}

Proof. Let $a \in W_J$ with $w \d(a) \in \cj_{I, \s}$ and $a \in
\cj_{J, \t}$. Then $\t(a)=a \i$ and $\s(w \d(a))=(w \d(a)) \i$.
Hence $\s(w)=\d(a) \i w \i a$. By Lemma 3.11, $a \t(\Phi_{I(J, w,
\d)})=\Phi_{I(J, w, \d)}$ Now write $a$ as $a=b u$ for $b \in
W_{I(J, w, \d)}$ and $u \in W_J$ with $u \t( \Phi_{I(J, w, \d)}^+)
\subset \Phi^+$. Then $u \t(\Phi_{I(J, w, \d)}^+)=\Phi_{I(J, w,
\d)}^+$.

By assumption, there exists $g \in G$ such that $g \i \s(g) \in B w
\d(a) \s(B)$. Notice that $$w \d(u) \s(\Phi_{I(J, w, \d)}^+)=w \d(u
\t( \Phi_{I(J, w, \d)}^+))=w \d(\Phi_{I(J, w, \d)}^+)=\Phi_{I(J, w,
\d)}^+.$$ Thus $w \s(a) \in W_{I(J, w, \d)} w \d(u)$ and $g \i \s(g)
\in P_{I(J, w, \d)} w \d(u) \s(P_{I(J, w, \d)})$. By Lemma 3.12,
there exists $l \in L_{I(J, w, \d)}$, such that $(g l) \i \s(g l)
\in B w \d(u) \s(B)$. In particular, $w \d(u) \in \cj_{I, \s}$.

Similarly, $u \in \cj_{J, \t}$. The existence is proved.

Now we prove the uniqueness. Assume that $u_1, u_2 \in W_J$ with
$u_1 \t(\Phi_{I(J, w, \d)}^+)=u_2 \t(\Phi_{I(J, w,
\d)}^+)=\Phi_{I(J, w, \d)}^+$ and $\s(w)=\d(u_1) \i w \i u_1=\d(u_2)
\i w \i u_2$. Set $v=u_1 u_2 \i$. Then $v w=w \d(v)$. By Lemma 3.10,
$v \in W_{I(J, w, \d)}$. Notice that $u_1=v u_2$ and $u_1
\t(\Phi_{I(J, w, \d)}^+)=u_2 \t(\Phi_{I(J, w, \d)}^+)=\Phi_{I(J, w,
\d)}^+$. Then $v=1$ and $u_1=u_2$. The corollary is proved. \qed

\subsection*{3.14} Unless otherwise stated, we fix $w \in
\cw(J, \s, \t)$ in the rest of this section. We will simply write
$I(J, w, \d)$ as $K$. Let $u$ be the unique element in $W_J$ with $u
\t( \Phi_K^+)=\Phi_K^+$, $w \d(u) \in \cj_{I, \s}$ and $u \in
\cj_{J, \t}$. Set
\begin{align*} G_w &=\{g \in G; g \i \s(g) \in P_K w {\d(u)} \s(P_K)\},
\\ L_w &=\{l \in L_J; l \i \t(l) \in (P_K \cap L_J) u \t(P_K \cap L_J)\}.\end{align*}

Let $v_1$ be a $G^{\s} \times P_K$-orbit on $G_w$ and $v_2$ be a $L_J^{\t}
\times (P_K \cap L_J)$-orbit on $L_w$.

Notice that $w \d(u) \s(\Phi_K^+)=\Phi_K^+$, By Lemma 3.12, there
exists $g \in v_1$ with $g \i \s(g) \in B w \d(u) B$. By Proposition
3.7, there exists $x_1 \in v_1$ with $x_1 \i \s(x_1) \in w \d(u) T$.
Similarly, there exists $x_2 \in v_2$ such that $x_2 \i \t(x_2) \in
u T$.

Set $\s'=\Ad(x_1 \i \s(x_1)) \circ \s$ and $\t'=\Ad(x_2 \i \t(x_2))
\circ \t$. Then $\s'$ is an involution on $G$, $\t'$ is an
involution on $L_J$ and $\s'(L_K)=\t'(L_K)=L_K$. Moreover, $x_1
G^{\s'} x_1 \i=G^\s$ and $x_2 L_K^{\t'} x_2 \i=L_K^\t$.

Set $H=P_K \cap \t'(P_K) \cap L_J$. Then $U_H=U_{P_K} \cap
\t'(U_{P_K}) \cap L_J$.

\begin{th15} Keep the notation of 3.14. Then $$\pi_J(G^{\s'} \cap P_K)
L_K L_J^{\t'}=H L_J^{\t'}.$$
\end{th15}

The proof will be given in 3.20. The key point is to show that
certain equations on $U$ have common solutions (see Lemma 3.19). We
use the exponential map $exp: \Lie(U) \rightarrow U$ to reduce this
problem to the problem of solving certain equations on $\Lie(U)$.
(Here we use the fact that $exp$ is an isomorphism when the
characteristic of $\mathbf k$ is $0$ or sufficiently large.)

The equations on $\Lie(U)$ that we need to solve are nonlinear
equations. We will use ``linear approximation'' to solve these
equations. We will provide the setting for ``linear approximation''
in Lemma 3.16 and we will prove the existence of common solutions
for the linear equations in Lemma 3.17.

\begin{th16} Set $$u_0=\Lie(U_H)=\bigl(\Lie(U_{P_K}) \cap \Lie(L_J)
\bigr) \cap \t' \bigl(\Lie(U_{P_K}) \cap \Lie(L_J) \bigr).$$ Define
$u_n=[u_0, u_{n-1}]$ for $n>0$ and $u'_n=u_n \oplus \Lie(U_{P_J})$
for $n \ge 0$. Then

(1) $u_{n+1} \subset u_n$ and $u'_{n+1} \subset u'_n$ for all $n \ge 0$.

(2) $[u'_0, u'_n] \subset u'_{n+1}$ for all $n \ge 0$.

(3) $u_n=\{0\}$ and $u'_n=\Lie(U_{P_J})$ for $n \gg 0$.

(4) $\t'(u_n)=u_n$ and $\s'(u_n) \subset u'_n$ for $n \ge 0$.
\end{th16}

Proof. We prove (1) by induction on $n \ge 0$. We have
$$[\Lie(U_{P_K}) \cap \Lie(L_J), \Lie(U_{P_{I(J, w, \d)}}) \cap
\Lie(L_J)] \subset \Lie(U_{P_K}) \cap \Lie(L_J).$$

Thus $u_1=[u_0, u_0] \subset \Lie(U_{P_K}) \cap \Lie(L_J)$.

Similarly, $u_1 \subset \t' \bigl(\Lie(U_{P_K}) \cap \Lie(L_J)
\bigr)$. Hence $u_1 \subset u_0$. Assume now that $n>0$ and $u_n
\subset u_{n-1}$. Then $u_{n+1}=[u_0, u_n] \subset [u_0,
u_{n-1}]=u_n$. From definition, we deduce that $u'_{n+1} \subset
u'_n$.

We prove (2). By definition, $$u'_{n+1} =u_{n+1} \oplus
\Lie(U_{P_J})=[u_0, u_n] \oplus \Lie(U_{P_J})=[u'_0, u'_n]+
\Lie(U_{P_J}).$$

(3) follows easily from definition.

We prove (4) by induction on $n \ge 0$.

By definition, \begin{gather*} u_0=\bigoplus_{\a, \t'(\a) \in
\Phi_J^+-\Phi_K^+} u_{\a}, \\ u'_0=\bigl(\bigoplus_{\a, \t'(\a) \in
\Phi_J^+-\Phi_K^+} u_{\a} \bigr) \bigoplus \bigl(\bigoplus_{\a \in
\Phi^+-\Phi_J^+} u_{\a} \bigr).
\end{gather*}

Hence $\t'(u_0)=u_0$. Let $\a$ be a root with $\a, \t'(\a) \in
\Phi_J^+-\Phi_K^+$, then $$\s'(\a)=w \d(u) \s(\a)=w \d \t'(\a) \in w
\d(\Phi_J^+-\Phi_K^+) \subset \Phi^+-\Phi_K^+.$$

If moreover, $\s'(\a) \in \Phi_J^+$, then $\d \i(w \i \a)=\t'
\s'(\a)\in \Phi_J$. Since $w \in W^{\d(J)}$ and $\a \in \Phi^+$, we
deduce that $\t' \s'(\a) \in \Phi_J^+$. Therefore $\s' (u_0) \subset
u'_0$.

Assume now that $n>0$ and that (4) holds when $n$ is replaced by
$n-1$. Then \begin{gather*} \t'(u_n)=\t' ([u_0, u_{n-1}])=[\t'(u_0),
\t'(u_{n-1})]=[u_0, u_{n-1}]=u_n, \\ \s'(u_n)=\s' ([u_0,
u_{n-1}])=[\s'(u_0), \s'(u_{n-1})] \subset [u'_0, u'_{n-1}] \subset
u'_n.\end{gather*} (4) is proved. \qed

\begin{th17} Keep the notation of Lemma 3.16.
Let $x \in L_K$ with $\s'(x)=x \i$. Set $\tilde \s'=\Ad(x) \circ
\s'$. Then for each $a \in u_n$ with $\t'(a)=-a$, there exists $b
\in u'_n$ such that $\tilde \s'(b)=b$ and
$\pi_J(b)-\t'(\pi_J(b))=a$.
\end{th17}

Proof. We show that

(a) For any $n \ge 0$, $\Ad(x) u'_n=u'_n$.

We use induction on $n$. Notice that $x \in L_K$ and $\t'(L_K)=L_K$.
Then
\begin{gather*} \Ad(x) \bigl(\Lie(U_{P_K}) \cap \Lie(L_J)\bigr)
=\Lie(U_{P_{I(J, w, \d)}}) \cap \Lie(L_J), \\ \Ad(x) \t'
\bigl(\Lie(U_{P_{I(J, w, \d)}}) \cap \Lie(L_J) \bigr)=\t'
\bigl(\Lie(U_{P_K}) \cap \Lie(L_J) \bigr), \\ \Ad(x)
\Lie(U_{P_J})=\Lie(U_{P_J}).
\end{gather*}
Thus $\Ad(x) u'_0=u'_0$. Assume now that $n>0$ and that (a) hold
when $n$ is replaced by $n-1$. Then $$\Ad(x) u'_n=\Ad(x) [u'_0,
u'_{n-1}]+\Ad(x) \Lie(U_{P_J})=[u'_0, u'_{n-1}]+
\Lie(U_{P_J})=u'_n.$$ (a) is proved.

(b) For any $n \ge 0$, $\tilde \s'(u_n) \subset u'_n$.

This follows from (a) and Lemma 3.16(4).

(c) Let $\a$ be a root. If $(w \d)^n \a \in \Phi_J^+$ for all $n \ge 0$, then $a \in \Phi_K^+$.

We have that $(w \d)^k \a=\a$ for some $k>0$. Then $w \d
\bigl(\sum_{i=1}^{k-1} (w \d)^i \a \bigr)=\sum_{i=1}^{k-1} (w \d)^i
\a$. Write $\sum_{i=1}^{k-1} (w \d)^i \a$ as $\sum_{j \in J} a_j
\a_j$. Then $$w \d \{j \in J; a_j \neq 0\}=\{j \in J; a_j \neq
0\}.$$ Hence if $a_j \neq 0$, then $j \in I(J, w, \d)$. Therefore
$\a \in \Phi^+_K$.

(c) is proved.

(d) $(\pi_J \tilde \s' \t')^n \Lie(U_{P_K})=\{0\}$ for $n \gg 0$.

Notice that $\pi_J \tilde \s' \t' \in \pi_J \Ad(L_K) \Ad(w g_D)$.
Thus $$(\pi_J \tilde \s' \t')^n \Lie(U_{P_K})=(\pi_J \Ad(w g_D))^n
\Lie(U_{P_K}).$$ Now (d) follows from (c).

Now set $$b=-\frac{1+\tilde \s'}{2} \t' \sum_{i \ge 0} (\pi_J \tilde \s' \t')^n a.$$

By (d), $b$ is well-defined. By (b) and Lemma 3.16(4), $b \in u'_n$.
Since $\tilde \s'$ is an involution, we have that $\tilde \s'(b)=b$.
Now \begin{align*} (1-\t') \pi_J(b) &=\sum_{i \ge 0} (\pi_J \tilde
\s' \t')^n a/2-\t' \sum_{i \ge 0} (\pi_J \tilde \s' \t')^n a/2  \\ &
\quad +\t' \sum_{i \ge 1} (\pi_J \tilde \s' \t')^n a/2-\sum_{i \ge
1} (\pi_J \tilde \s' \t')^n a/2 \\ & =a/2-\t'(a/2)=a
\end{align*} \qed

\subsection*{3.18} Let us recall the Campell-Hausdorff formula.

Let $exp: \Lie(U) \rightarrow U$ be the exponential map. Then for
$X, Y \in \Lie(U)$, $exp(X)exp(Y)=exp(X+Y+\sum_{n>1} f_n(X, Y))$,
where $$f_n(X, Y)=\sum\limits_{\scriptstyle r_i+t_i>0; \atop
\scriptstyle r_1+t_1+\cdots+r_s+t_s=n-1} a_{r_1, t_1, \cdots, r_s,
t_s} (ad X)^{r_1} (ad Y)^{t_1} \cdots (ad X)^{r_s} (ad Y)^{t_s} Y$$
and the coefficients $a_{r_1, t_1, \cdots, r_s, t_s}$ only depends
on $r_1, t_1, \cdots, r_s, t_s$.

\begin{th18} Keep the notation of Lemma 3.16 and
3.17. Let $a \in u_0$ with $\t'(a)=-a$. Then there exists $b \in
u'_0$ such that $\tilde \s'(b)=b$ and $exp(\pi_J(b))
exp(-\t'(\pi_J(b)))=exp(a)$.
\end{th18}

Proof. For $b \in \Lie(P_J)$, we simply write $\pi_J(b)$ as $\bar
b$. It suffices to prove the following statement: for each $n$,
there exists $b_n \in u'_0$, such that $$\tilde \s'(b_n)=b_n \text{
and } exp(\bar b_n) exp(-\t'(\bar b_n)) \in exp (a+u_n).$$

We prove by induction on $n$. For $n=0$, we may choose $b_0=0$.
Assume that $n>0$ and that $exp(\bar b_{n-1}) exp(-\t'(\bar
b_{n-1})) \in exp (a+u_{n-1})$ for $\tilde \s'(b_{n-1})=b_{n-1} \in
u'_0$. By Campell-Hausdorff formula, $$\bar b_{n-1}-\t'(\bar
b_{n-1})+\sum_{i=2}^{n-1} f_i(\bar b_{n-1}, -\t'(\bar b_{n-1})) \in
a+u_{n-1}.$$ Thus $\bar b_{n-1}-\t'(\bar b_{n-1})+\sum_{i=2}^n
f_i(\bar b_{n-1}, -\t'(\bar b_{n-1})) \in a+a_{n-1}+u_n$ for some
$a_{n-1} \in u_{n-1}$.

Since $\t' \bigl(exp(\bar b) exp(-\t'(\bar b)) \bigr)=\bigl(exp(\bar
b) exp(-\t'(\bar b)) \bigr) \i$, we have that
$\t'(a_{n-1})=-a_{n-1}$. By the definition of $f_i$ and Lemma
3.16(4),
\begin{align*} f_i(\bar b_{n-1}+u_{n-1}, -\t'(\bar b_{n-1}+u_{n-1}))
&=f_i(\bar b_{n-1}+u_{n-1}, -\t'(\bar b_{n-1})+u_{n-1}) \\
& \subset f_i(\bar b_{n-1}, -\t'(\bar b_{n-1}))+[u_0, u_{n-1}] \\
&=f_i(\bar b_{n-1}, -\t'(\bar b_{n-1}))+u_n. \end{align*}

By Lemma 3.17, there exists $b'_{n-1} \in u'_{n-1}$, such that
$\tilde \s'(b'_{n-1})=b'_{n-1}$ and $b'_{n-1}-\t'(b'_{n-1}) \in
a_{n-1}+\Lie(U_{P_J})$. Set $b_n=b_{n-1}-b'_{n-1}$. Then $$\bar
b_n-\t'(\bar b_n)+\sum_{i=2}^n f_i(\bar b_n, -\t'(\bar b_n)) \in
a+u_n.$$

In other words, $exp(\bar b_n) exp(-\t'(\bar b_n)) \in exp (a+u_n)$.
\qed

\subsection*{3.20} The proof of Proposition 3.15. By
definition $$\pi_J (\Lie(U_{P_K}) \cap \s'\Lie(U_{P_K})) \subset
\bigoplus_{\a, \s'(\a) \in \Phi^+_J-\Phi^+_K} u_{\a}.$$

Let $\a$ be a root with $\a, \s'(\a) \in \Phi^+_J-\Phi^+_K$. It is
easy to see that $\t'(\a) \in \Phi_J \cap \t' \s'(\Phi^+_J)=\Phi_J
\cap \d \i w \i(\Phi^+_J)$. Since $w \in W^{\d(J)}$, we deduce that
$\t'(\a) \in \Phi^+_J$. Notice that $\t'(\Phi_K)=\Phi_K$. Thus
$\t'(\a) \in \Phi^+_J-\Phi^+_K$. Hence $\t' \pi_J (\Lie(U_{P_K})
\cap \s'\Lie(U_{P_K})) \subset \Lie(U_{P_K})$. Therefore $$\t'
\pi_J(P_K \cap \s'(P_K)) =\t'(L_K) \t' \pi_J(P_K \cap \s'(P_K))
\subset P_K.$$ In other words, $\pi_J(P_K \cap \s'(P_K)) \subset
\t'(P_K)$. Therefore $$\pi_J(P_K \cap \s'(P_K)) \subset H.$$ In
particular, $\pi_J(G^{\s'} \cap P_K) L_K L_J^{\t'} \subset H
L_J^{\t'}$.

Let $l \in L_K$ and $p \in U_H$. Then $p \t'(p) \i=\exp(a)$ for some
$a \in \Lie(U_H)$ with $\t'(a)=-a$. Set $\tilde \s'=\Ad(l \s'(l) \i)
\circ \s'$. By Lemma 3.19, there exists $g \in G^{\tilde \s'} \cap
U_{P_K}$ such that $\pi_J(g) \t'(\pi_J(g)) \i=\exp(a)$. So $p \in
\pi_J(G^{\tilde \s'} \cap U_{P_K}) L_J^{\t'}$. Thus $l p \in l
\pi_J(l \i G^{\s'} l \cap U_{P_K}) L_J^{\t'}=l \pi_J(G^{\s'} \cap
U_{P_K}) L_J^{\t'}$. The proposition is proved.

\subsection*{3.21} Keep the notation of 3.14. Set $$X_{J, \t; v_1, v_2}
=x_1 G^{\s'} L_K U_{P_J} L_J^{\t'} x_2 \i /U_{P_J} L_J^\t=G^\s x_1
L_K x_2 \i U_{P_J} L_J^\t/U_{P_J} L_J^\t.$$

We call $X_{J, \t; v_1, v_2}$ a $G^\s$-stable piece in $X_{J, \t}$.

\begin{th19} The variety $X_{J, \t; v_1, v_2}$ consists of the element $g l \i U_{P_J} L_J^\t$,
where $g \in v_1$ with $g \i \s(g) \in w \d(u) T$ and $l \in v_2$
with $l \i \t(l) \in (P_K \cap L_J) u$. In particular, $X_{J, \t;
v_1, v_2}$ is independent of the choice of $x_1$ and $x_2$.
\end{th19}

Proof. It is easy to see that $g \i \s(g) \in w \d(u) T$ for $g
\in x_1 G^{\s'}=G^\s x_1$ and $l \i \t(l) \in L_K u$ for $l \in x_2 L_K$.

Let $g=k x_1 p$ for $k \in G^{\s}$ and $p \in P_K$. If moreover, $g
\i \s(g) \in w \d(u) T$, then $\s'(p) \in P_K$. By 3.20, $\pi_J(p)
\in H$ and $p \in H U_{P_J}$.

Let $l=k' x_2 p'$ for $k' \in L_J^\t$ and $p' \in P_K \cap L_J$. If
moreover, $l \i \t(l) \in (P_K \cap L_J) u$, then $(p') \i \t'(p')
\in P_K \cap L_J$ and $\t'(p') \in P_K \cap L_J$. Thus $p' \in H$.

By Proposition 3.15, $g l \i \in G^\s x_1 H U_{P_J} H x_2 \i
L_J^\t=G^\s x_1 L_K x_2 \i U_{P_J} L_J^\t$. Hence $g l \i U_{P_J}
L_J^\t \in X_{J, \t; v_1, v_2}$. \qed

\begin{th20} We have that $$X_{J, \t}=\bigsqcup_{w \in \cw(J,
\s, \t)} \bigsqcup_{v_1 \in G^{\s} \setminus G_w/P_K} \bigsqcup_{v_2
\in L_J^{\t} \setminus L_w/(P_K \cap L_J)} X_{J, \t; v_1, v_2}.$$
\end{th20}

Proof. By Proposition 3.9, $X_{J, \t}=\bigsqcup_{w \in \cw(J, \s,
\t)} (X_{J, \t} \cap Z_{J, y, D; w})$.

Let $w \in \cw(J, \s, \t)$, $v_1 \in G^{\s} \setminus G_w/P_K$ and
$v_2 \in \in L_J^{\t} \setminus L_w/(P_K \cap L_J)$. Then
\begin{align*} X_{J, \t; v_1, v_2} & \subset (G^\s)_{\D} (\s(x_1 L_K x_2 \i),
x_1 L_K x_2 \i) \cdot h_{J, y, D} \\ &=(G^\s)_{\D} (x_1 w \d(u), x_1
L_K u L_K) \cdot h_{J, y, D} \\ &=(G^\s)_{\D} (x_1 w, x_1 L_K) \cdot
h_{J, y, D} \subset Z_{J, y, D; w}. \end{align*}

Let $z \in Z_{J, y, D; w} \cap X_{J, \t}$. By the proof of
Proposition 3.9, $z$ can be written as $z=(g b w, g) \cdot h_{J,
D}$, where $b \in B$, $g \in G$ with $g \i \s(g) \in B w \d(u') B$
for some $u' \in \cj_{J, \t}$. By the proof of Corollary 3.13, $u'
\in W_K u$. Thus $g \i \s(g) \in P_K w \d(u) \s(P_K)$. In other
words, $g \in v_1$ for some $v_1 \in G^{\s} \setminus G_w/P_K$. Now
\begin{align*} z & \in \bigl((G^{\s})_{\D} (x_1 P_K w, x_1 P_K) \cdot h_{J, y, D}\bigr) \cap X_{J, \t} \\
&=(G^{\s})_{\D} \bigl((\s(x_1) {\d(u)} \i w \i P_K w, x_1 P_K) \cdot
h_{J, y, D} \cap X_{J, \t} \bigr).
\end{align*}

Similar to the proof of Proposition 3.9, $z=(\s(k x_1 l \i), k x_1 l
\i) \cdot h_{J, y, D}$, where $k \in G^\s$ and $l \in v_2$ with $l
\i \t(l) \in (P_K \cap L_J) u$. Thus $l \in v_2$ for some $v_2 \in
L_J^{\t} \setminus L_w/(P_K \cap L_J)$. By Lemma 3.22, $z \in X_{J,
t; v_1, v_2}$.

Hence $X_{J, \t}=\bigsqcup_{w \in \cw(J, \s, \t)} \bigcup_{v_1 \in
G^{\s} \setminus G_w/P_K} \bigcup_{v_2 \in L_J^{\t} \setminus
L_w/(P_K \cap L_J)} X_{J, \t; v_1, v_2}$.

Let $w \in \cw(J, \s, \t)$, $v_1, v_1' \in G^{\s} \setminus G_w/P_K$
and $v_2, v_2' \in \in L_J^{\t} \setminus L_w/(P_K \cap L_J)$.
Assume that $X_{J, \t; v_1, v_2} \cap X_{J, \t; v_1', v_2'} \neq
\varnothing$, i.e., there exists $g_1 \in v_1, g_2 \in v_1'$, $l_1
\in v_2, l_2 \in v_2'$ such that $g_1 \i \s(g_1), g_2 \i \s(g_2) \in
w \d(u) T$, $l_1 \i \t(l_1), l_2 \i \t(l_2) \in (P_K \cap L_J) u$
and $$(\s(g_1 l_1 \i), g_1 l_1 \i) \cdot h_{J, y, D}=(\s(g_2 l_2
\i), g_2 l_2 \i) \cdot h_{J, y, D}.$$

Notice that
\begin{gather*} (\s(g_1 l_1 \i), g_1 l_1 \i) \cdot h_{J, y, D} \in
(g_1 w, g_1 P_K) \cdot h_{J, y, D}; \\ (\s(g_2 l_2 \i), g_2 l_2 \i)
\cdot h_{J, y, D} \in (g_2 w, g_2 P_K) \cdot h_{J, y, D}.
\end{gather*}

By 1.3 (2), $g_1 \in g_2 P_K$. In particular, $v_1=v_1'$. By the
proof of Lemma 3.22, $g_1, g_2 \in G^\s x_1 H U_{P_J}$. Thus
$$g_1 \i g_2 \in H U_{P_J} G^{\s'} H U_{P_J} \cap P_K =H U_{P_J}
(G^{\s'} \cap P_K) H U_{P_J}=H U_{P_J}.$$

Now \begin{align*} (1, l_1 \i \t(l_1)) \cdot h_{J, y, D} &=(1,
\pi_J(g_1 \i g_2) l_2 \i \t(l_2) \t(\pi_J(g_1 \i g_2)) \i) \cdot
h_{J, y, D} \\ &=(1, (l_2 \pi_J(g_2 \i g_1)) \i \t(l_2 \pi_J(g_2 \i
g_1))) \cdot h_{J, y, D}. \end{align*}

Notice that $l_1, l_2 \pi_J(g_2 \i g_1) \in L_J$. Then $$l_1 \i
\t(l_1)=(l_2 \pi_J(g_2 \i g_1)) \i \t(l_2 \pi_J(g_2 \i g_1)).$$
Hence $l_1 \in l_2 \pi_J(g_2 \i g_1) L_J^\t$ and $v_2=v_2'$. \qed

\subsection*{3.24} By the proof of Theorem 3.23, each $G^\s$-stable piece
in $X_{J, \t}$ is an irreducible component of the intersection of
$X_{J, \t}$ with some $G$-stable piece in $Z_{J, y, D}$. In
particular, in the Example 3.6 (1), the $\mathbf G^\s$-stable pieces
in $\mathbf X_{\mathbf J, \s}$ are just the $G$-stable pieces in
$Z_{J, 1, G}$ (see remark of Proposition 3.9).

\subsection*{3.25} For any variety $X$ with the action of $G^{\s'} \cap P_K$ (resp. $L_K^{\s'}$), we denote by $ind^1(X)$ (resp.
$ind^2(X)$) the quotient of $G^{\s'} \times X$ modulo the action of
$G^{\s'} \cap P_K$ (resp. $L_K^{\s'}$) defined by $g \cdot (g',
x)=(g' g \i, g \cdot x)$.

Recall that $\pi_J(G^{\s'} \cap P_K) \subset H$ (see 3.20). Unless
otherwise stated, the action of $G^{\s'} \cap P_K$ on $H L_J^{\t'}$
is defined by $g \cdot h=\pi_J(g) h$ and the action of $G^{\s'} \cap
P_K$ on $L_K$ is defined by $g \cdot l=\pi_K(g) l$.

\begin{th21} The map $G^{\s'} \times H L_J^\t \rightarrow G$ defined by $(g, l) \mapsto x_1 g
l x_2 \i$ induces an isomorphism $\kappa: ind^1(H
L_J^{\t'}/L_J^{\t'}) \xrightarrow \simeq X_{J, \t; v_1, v_2}$.
\end{th21}

Proof. By 1.3 (2), the map $\pi: Z_{J, y, D; w} \rightarrow G/P_K$
defined by $$\pi \bigl((g p  w, g) \cdot h_{J, y, D} \bigr)=g P_K$$
for $g \in G$ and $p \in P_K$ is well-defined and is a
$G$-equivariant morphism. Then its restriction to $(x_1 \i, x_1 \i)
X_{J, \t; v_1, v_2}$ is a $G^{\s'}$-equivariant morphism.

Notice that $(x_1 \i, x_1 \i) X_{J, \t; v_1, v_2} \subset
(G^{\s'})_{\D} (P_K w, P_K) \cdot h_{J, y, D}$. Then $$\pi((x_1 \i,
x_1 \i) X_{J, \t; v_1, v_2}) \subset G^{\s'} P_K/P_K \cong
G^{\s'}/{G^{\s'} \cap P_K}.$$

Since $G^{\s'}$ acts transitively on $G^{\s'}/{G^{\s'} \cap P_K}$, we have that
$$\pi((x_1 \i, x_1 \i) X_{J, \t; v_1, v_2})=G^{\s'}/{G^{\s'} \cap P_K}.$$

By the proof of Theorem 3.23,
\begin{align*} & (\pi \mid_{(x_1 \i, x_1 \i) X_{J, \t; v_1, v_2}}) \i
(G^{\s'} \cap P_K)=(P_K w, P_K) \cdot h_{J, y, D} \cap (x_1 \i, x_1
\i) X_{J, \t; v_1, v_2} \\ &=\{(w, l \i \t(l)) \cdot h_{J, y, D}; l
\in v_2 \text{ with } l \i \t(l) \in (P_K \cap L_J) u\}.
\end{align*}

We have shown in the proof of Lemma 3.22 that $l \in v_2$ with $l \i
\t(l) \in (P_K \cap L_J) u$ if and only if $l \in x_2 \i L_J^{\t'}
H$. Thus $$(\pi \mid_{(x_1 \i, x_1 \i) X_{J, \t; v_1, v_2}}) \i
(G^{\s'} \cap P_K) \cong H L_J^{\t'}/L_J^{\t'}.$$ Now the lemma
follows from Lemma 3.3. \qed

\begin{th22} Keep the notation of 3.14. Then

(1) $U_H/U_H^{\t'}$ is an affine space.

(2) The map $L_K \times U_H \rightarrow H$ defined by $(l, p)
\mapsto l p$ induces an isomorphism $L_K/L_K^{\t'} \times
U_H/U_H^{\t'} \rightarrow H L_J^{\t'}/L_J^{\t'}$.
\end{th22}

Proof. (1) Notice that $\t'$ is an order-2 linear endomorphism on
$\Lie(U_H)$. So $\Lie(U_H)=\Lie(U_H)^{\t'} \oplus \Lie(U_H)^{1-\t'}$
and $\Lie(U_H)/\Lie(U_H)^{\t'}$ is an affine space.

Moreover, the isomorphism $\exp: \Lie(U_H) \rightarrow U_H$ induces
an isomorphism $\Lie(U_H)/\Lie(U_H)^{\t'} \rightarrow
U_H/U_H^{\t'}$. Part (1) is proved.

Part (2) is obvious. \qed

\

The following result is an easy consequence of the above lemma and
Proposition 3.15.

\begin{th23} The projection map $H L_J^{\t'}/L_J^{\t'} \rightarrow
L_K/L_K^{\t'}$ induces an affine space bundle map $$\v: ind^1(H
L_J^{\t'}/L_J^{\t'}) \rightarrow ind^1(L_K/L_K^{\t'}).$$ Moreover,
this map induces a bijection from the set of $G^{\s'}$-orbits on
$ind^1(H L_J^{\t'}/L_J^{\t'})$ to the set of $G^{\s'}$-orbits on
$ind^1(L_K/L_K^{\t'})$.
\end{th23}

\

Notice that the map $\kappa$ induces a bijection between the set of
$G^\s$-orbits on $X_{J, \t; v_1, v_2}$ and the set of
$G^{\s'}$-orbits on $ind^1(H L_J^{\t'}/L_J^{\t'})$. Then

\begin{th24} There is a bijection between the set of $G^\s$-orbits on
$X_{J, \t; v_1, v_2}$ and the set of $L_K^{\s'}$-orbits on
$L_K/L_K^{\t'}$.
\end{th24}

\subsection*{3.30} In the rest of this section, we assume that $G$ is adjoint.
We assume furthermore that for $\a \in \Phi^+$, either $\t(\a)=\a$
or $\t(\a) \in -\Phi^+$. Let $I$ be the set of simple roots and
$I_0$ be the set of simple roots that are fixed by $\t$.

Let $\overline{G/G^{\t}}$ be the De Concini-Procesi compactification
of $G/G^{\t}$. Then $\overline{G/G^{\t}}$ is a smooth, projective
variety that contains $G/G^{\t}$ as an open subvariety. The
$G^\s$-action on $G/G^{\t}$ extends in a unique way to a
$G^\s$-action on $\overline{G/G^{\t}}$. Moreover,
$$\overline{G/G^{\t}}=\bigsqcup_{I_0 \subset J \subset I;
\t(\Phi_J)=\Phi_J} \bar X_{J, \t}.$$ Here $\bar X_{J, \t}$ is the
quotient space $G \times_{P_J} (G_J/G_J^{\t})$, where
$G_J=L_J/Z^0(L_J)$ and $P_J$ acts on $G$ on the right and acts on
$G_J/G_J^{\t}$ via the quotient $L_J$ (see \cite[1.4]{Sp1}).

Let $p_J: X_{J, \t} \rightarrow \bar X_{J, \t}$ be the projection
map. Set $\bar X_{J, \t; v_1, v_2}=p_J(X_{J, \t; v_1, v_2})$. We
call $\bar X_{J, \t; v_1, v_2}$ a $G^\s$-stable piece in
$\overline{G/G^\t}$.

\begin{th241} We have that
$$\overline{G/G^{\t}}=\bigsqcup_{I_0 \subset J \subset I; \t(\Phi_J)=\Phi_J}
\bigsqcup_{w \in \cw(J, \s, \t)} \bigsqcup\limits_{\scriptstyle v_1
\in G^{\s} \setminus G_w/P_{I(J, w, \d)}; \atop \scriptstyle v_2 \in
L_J^{\t} \setminus L_w/(P_{I(J, w, \d)} \cap L_J)} \bar X_{J, \s;
v_1, v_2}.$$

Let $I_0 \subset J \subset I$ with $\t(\Phi_J)=\Phi_J$, $w \in
\cw(J, \s, \t)$, $v_1 \in G^{\s} \setminus G_w/P_{I(J, w, \d)}$ and
$v_2 \in L_J^{\t} \setminus L_w/(P_{I(J, w, \d)} \cap L_J)$. We use the same notation as in 3.14. Then
$\bar X_{J, \t; v_1, v_2}$ is an affine space bundle over $G^{\s'}
\times_{G^{\s'} \cap P_K} L_K/L_K^{\t'} Z^0(L_J)$. Moreover, this map induces
a bijection between the set of $G^{\s}$-orbits on
$\bar X_{J, \t; v_1, v_2}$ and the set of $\bigl(L_K/Z^0(L_J) \bigr)^{\s'}$-orbits on
$\bigl(L_K/Z^0(L_J) \bigr)/\bigl(L_K/Z^0(L_J) \bigr)^{\t'}$.
\end{th241}

\section{The character sheaves}

\subsection{} We follow the notation of \cite{BBD} and \cite{BL}. Let $X$ be an
algebraic variety over $\mathbf k$ and $l$ be a fixed prime number
invertible in $\mathbf k$. We write $\cd(X)$ instead of $\cd^b_c(X,
\bar{\mathbb Q}_l)$. If $f: X \rightarrow Y$ is a smooth morphism
with connected fibres of dimension $d$, then we set $\tilde
f(C)=f^*(C)[d]$ for any perverse sheaf $C$ on $Y$.

Let $K$ be an algebraic group defined over $\mathbf k$. If $K$ acts
on $X$, we denote by $\cd_K(X)$ the equivariant derived category of
$X$.

\subsection{} Let $C_i \in \cd_K(X)$ for $i=1, 2, \cdots, n$. For $C \in \cd_K(X)$,
we write $C \in <C_i; i=1, 2, \cdots, n>$ if there exist $m>n$ and
$C_{n+1}, \cdots, C_m \in \cd_K(X)$ such that $C_m=C$ and for each
$n+1 \le i \le m$, there exists $1 \le j, k<i$ such that $(C_j, C_i,
C_k)$ is a distinguished triangle in $\cd_K(X)$.

If $X, C, C_i$ are as above and $Y \xrightarrow{f} X \xrightarrow{g}
Z$ are $K$-equivariant morphisms, then \begin{align*} & f^*(C) \in
<f^*(C_i); i=1, 2, \cdots, n>, \tag{a} \\ & g_!(C) \in <g_!(C_i);
i=1, 2, \cdots, n>. \tag{b} \end{align*}

If $X=\sqcup_{1 \le i \le n} X_i$ is a partition of $X$ into locally
closed $K$-stable subvarieties such that $\sqcup_{1 \le i \le k}
X_i$ is closed in $X$ for any $1 \le k \le n$. We denote by $j_i:
X_i \rightarrow X$ the inclusion maps. Let $C \in \cd_K(X)$. Then
\begin{align*} C \in <(j_i)_! (j_i)^* (C); i=1, 2, \cdots, n>.
\tag{c}\end{align*}

In the case when $K$ is a trivial group, the notation above is
slightly different from the one defined in \cite[32.15]{L2}. Namely,
$C \in <C_i; i=1, 2, \cdots, n>$ if and only if there exists a
sequence $\{C'_i; i \in \mathbf Z\}$ of objects in $\cd(X)$ such
that $C'_i \in \{0, C_1, C_2, \cdots, C_n\}$ for all $i$, $C'_i=0$
for all but finitely many $i$ and $C \Bumpeq \{C'_i; i \in \mathbf
Z\}$.

\subsection{} Let $T$ be a torus. Let $\ck(T)$ be the set of
isomorphism classes of Kummer local systems on $T$, i. e., the set
of isomorphism classes of $\bar{\mathbb Q}_l$-local systems $\cl$ of
rank one on $T$, such that $\cl^{\otimes m} \cong \bar{\mathbb Q}_l$
for some integer $m \ge 1$ invertible in $\mathbf k$.

Let $X$ be a variety with free $T$-action $a: T \times X \rightarrow
X$. For $\cl \in \ck(T)$, we denote by $\cd^{\cl}(X)$ the full
subcategory of $\cd(X)$ with objects $A \in \cd(X)$ such that $a^* A
\cong \cl \boxtimes A$. If moreover, we have an action of an
algebraic group $K$ on $X$ that commutes with the action of $T$, we
denote by $\cd^{\cl}_K(X)$ the full subcategory of $\cd_K(X)$ with
objects $A \in \cd(X)$ such that the image of $A$ in $\cd(X)$ is in
$\cd^{\cl}(X)$.

\subsection{} Let $J \subset I$. Denote by $Y_J$ the quotient
of $G/U \times L_J/U_J$ modulo the diagonal $T$-action on the right.
Consider the diagram \[\xymatrix{G/U_{P_J} & G/U_{P_J} \times
L_J/B_J \ar[l]_(.6){p_J} \ar[r]^(.7){q_J} & Y_J}\] where $p_J$ is
the projection and $q_J(g, l)=(g l U, l U_J) T$. Then $p_J$ is
proper and $q_J$ is a smooth morphism with fibres isomorphic to
$U_J$.

Define the $G^\s \times L_J^\t$-action on $G/U_{P_J}$ by $(g, l)
\cdot g'=g g' l \i$, on $G/U_{P_J} \times L_J/U_J$ by $(g, l) \cdot
(g', l')=(g g' l \i, l l')$ and on $Y_J$ by $(g, l) \cdot (g' U, l'
U_J)T=(g g' U, l l' U_J)T$. Then the maps $p_J$ and $q_J$ are $G^\s
\times L_J^\t$-equivariant.

Define $CH_J=(p_J)_! (q_J)^*: \cd_{G^\s \times L_J^{\t}}(Y_J)
\rightarrow \cd_{G^\s \times L_J^\t}(G/U_{P_J})$ and $HC_J=(q_J)_*
(p_J)^!: \cd_{G^\s \times L_J^\t}(G/U_{P_J}) \rightarrow \cd_{G^\s
\times L_J^\t}(Y_J)$.

It is easy to see that the functor $CH_J$ is the left adjoint of
$HC_J$.

In the case when $J=I$, $CH_I$ is just the character functor defined
in \cite[8.4]{Gi} and $HC_I$ is just the Harish-Chandra functor
defined in {\it loc.cit}.

We will call $CH_J$ a (parabolic) character functor and $HC_J$ a (parabolic) Harish-Chandra functor.

\begin{th25} Let $A \in \cd_{G^\s \times L_J^\t}(G/U_{P_J})$. Then $A$ is a direct
summand of $CH_J \circ HC_J(A)$.
\end{th25}

\begin{rmk} The argument is inspired by \cite[8.5.1]{Gi} and \cite[3.6]{MV}.
\end{rmk}

Proof. Set $Z=\{(u, y); u \in L_J, y \in L_J/B_J, u \in {}^y U_J\}$.
The second projection $pr_2: Z \rightarrow L_J/B_J$ is a fibration
with fibres $U_J$. The first projection $pr_1: Z \rightarrow L_J$ is
the ``Springer resolution'' of the unipotent variety of $L_J$.

Consider the following commuting diagram

\[\xymatrix {& G/U_{P_J} \times L_J/B_J \ar[dl]_{p_J} \ar[d]_{q_J} & G/U_{P_J}
\times Z \ar[l]_(.4){id \times pr_2}  \ar[r]^{id \times pr_1}
\ar[d]_h & (G/U_{P_J}) \times L_J \ar[d]^m \\ G/U_{P_J} & Y_J &
G/U_{P_J} \times L_J/B_J \ar[l]_(.6){q_J} \ar[r]^(.6){p_J} &
G/U_{P_J}}\] where $h(g, u, y)=(gu, y)$ and $m(g, l)=g l$.

It is easy to see that the square $(h', q_J, h, q_J)$ is Cartesian.
Now
\begin{align*} CH_J \circ HC_J(A) &=(p_J)_* (q_J)^! (q_J)_* (p_J)^!(A)[-2 d] (-d) \\ &
=(p_J)_* h_* (id \times pr_2)^! (p_J)^!(A)[-2 d](-d) \\ &=m_* (id
\times pr_1)_* \bigl(A \boxtimes \bar{\mathbb Q}_l[2 d](d) \bigr) \\
&=m_* \bigl(A \boxtimes (pr_1)_* \bar{\mathbb Q}_l [2 d](d) \bigr)
\end{align*} Here $d=\dim(U_J)$ and $(-d)$, $(d)$ are Tate twists.

It is known that $(pr_1)_* \bar{\mathbb Q}_l [2 d](d)$ is a
semisimple perverse sheaf on $L_J$ and the skyscraper sheaf
$(\bar{\mathbb Q}_l)_e$ at the identity point of $L_J$ is a direct
summand of $(pr_1)_* \bar{\mathbb Q}_l [2 \dim(U_J)]$. Hence
$A=m_*(A \boxtimes (\bar{\mathbb Q}_l)_e)$ is a direct summand of
$CH_J \circ HC_J(A)$. \qed

\begin{th26} We define the action of $T$ on $Y_J$
by $t (g U, l U_J)T=(g t U, l U_J)T$. Let $\cl \in \ck(T)$ and $A
\in \cd_{G^\s \times L_J^\t}^{\cl}(Y_J)$. Then
$$HC_J \circ CH_J(A) \in \sum_{w \in W_J} \cd_{G^\s \times L_J^\t}^{w^* \cl}(Y_J).$$
\end{th26}

\begin{rmk} This result is inspired by \cite[Proposition 1.2]{Gr}.
\end{rmk}

Proof. Set $Z=G/U_{P_J} \times L_J/B_J \times L_J/B_J$. Consider the
following diagram

\[
\xymatrix{ & & Z \ar[dl]_{a} \ar[dr]^{b} & &
\\ & G/U_{P_J} \times L_J/B_J \ar[dl]_{q_J} \ar[dr]^{p_J} & & G/U_{P_J} \times L_J/B_J
\ar[dl]_{p_J} \ar[dr]^{q_J} &
\\ Y_J & & G/U_{P_J} & & Y_J}
\] where $a(g, l, l')=(g, l)$ and $b(g, l, l')=(g, l')$.

It is easy to see that the square $(a, p_J, b, p_J)$ is Cartesian.
Notice that $p_J$ is proper, then
\begin{align*} HC_J \circ CH_J(A) &=(q_J)_* (p_J)^! (p_J)_!
(q_J)^*=(q_J)_* a_* b^! (q_J)^!(A) \\ &=(q_J \circ a)_* (q_J \circ
b)^!(A). \end{align*}

Now we have a partition $Z=\sqcup_{w \in W_J} Z_w$, where $$Z_w=\{(g
U_{P_J}, l_1 B_J, l_2 B_J); l_1, l_2 \in L_J, l_1 \i l_2 \in B_J w
B_J\}.$$ By 4.2(c), $HC_J \circ CH_J(A) \in <((q_J \circ
a)\mid_{Z_w})_* ((q_J \circ b)\mid_{Z_w})^! (A); w \in W_J>$.

Set \begin{align*} Z'_w=\{\bigl( & (x U, y U_J)T, (a U, b U_J) T
\bigr) \in Y_J \times Y_J; \\ & U_J x \i a U_J=U_J y \i b U_J
\subset B_J w B_J\}. \end{align*} Define the map $\pi_w: Z_w
\rightarrow Z'_w$ by $$(g U_{P_J}, l_1 B_J, l_2 B_J) \mapsto \bigl(
(g l_1 U, l_1 U_J)T, (g l_2 U, l_2 U_J) T \bigr).$$ Then it is easy
to see that $\pi_w$ is an affine space bundle map with fibres
isomorphic to $U_J \cap {}^w U_J$.

Let $p_w: Z'_w \rightarrow Y_J$ be the projection to the first
factor and $p'_w: Z'_w \rightarrow Y_J$ be the projection to the
second factor. Then $(q_J \circ a)\mid_{Z_w}=p_w \circ \pi_w$ and
$(q_J \circ b)\mid_{Z_w}=p'_w \circ \pi_w$. Now \begin{align*} ((q_J
\circ a)\mid_{Z_w})_* ((q_J \circ b)\mid_{Z_w})^! (A) &=(p_w)_*
(\pi_w)_* (\pi_w)^! (p'_w)^! (A) \\ &=(p_w)_* (p'_w)^!(A)[2 d](d).
\end{align*} Here $d=\dim(U_J \cap {}^w U_J)$ and $(d)$ is Tate twist.

Define the $T$-action on $Z'_w$ by $$t \cdot \bigl( (x U, y U_J)T,
(a U, b U_J) T \bigr)=\bigl( (x w t w \i U, y U_J)T, (a t U, b U_J)
T \bigr).$$ Then $p'_w$ is $T$-equivariant and $p_w$ is
$T$-equivariant with respect to the twisted $T$-action on $Y_J$
defined by $t \cdot (x U, y U_J) T=(x w t w \i U, y U_J) T$. Thus
$(p_w)_* (p'_w)^!(A) \in \cd_{G^\s \times L_J^\s}^{w \i \cl}(Y_J)$.
\qed

\subsection*{4.7} Let $\cd^{cs}_{G^\s \times L_J^\t}(Y_J)$ be the
full subcategory of $\cd_{G^\s \times L_J^\t}(Y_J)$ with objects in
$\oplus_{\cl \in \ck(T)} \cd^{\cl}_{G^\s \times L_J^\t}(Y_J)$. By
\cite[5.3]{BL}, $CH_J(A)$ is semisimple for simple perverse sheaf $A
\in \cd^{cs}_{G^\s \times L_J^\t}(Y_J)$. Let $\cc_{G^\s \times
L_J^\t}(G/U_{P_J})$ be the set of (isomorphism classes) of simple
perverse sheaves that are a constituent of $CH_J(A)$ for some $A \in
\cc_{G^\s \times L_J^\t}(Y_J)$. The elements in $\cc_{G^\s \times
L_J^\t}(G/U_{P_J})$ are called (parabolic) character sheaves on
$G/U_{P_J}$. Let $\cd^{cs}_{G^\s \times L_J^\t}(G/U_{P_J})$ be the
full subcategory of $\cd_{G^\s \times L_J^\t}(G/U_{P_J})$ consisting
of objects whose perverse constituents are contained in $\cc_{G^\s
\times L_J^\t}(G/U_{P_J})$. Then it is easy to see that
$CH_J(\cd_{G^\s \times L_J^\t}^{cs}(Y_J)) \subset \cd_{G^\s \times
L_J^\t}^{cs}(G/U_{P_J})$.

Similar to \cite[Proposition 6.7(a)]{L3}, we have the following
result.

\begin{th27} Let $A \in \cd_{G^\s \times L_J^\t}(G/U_{P_J})$. Then
$A \in \cd_{G^\s \times L_J^\t}^{cs}(G/U_{P_J})$ if and only if
$HC_J(A) \in \cd_{G^\s \times L_J^\t}^{cs}(Y_J)$.
\end{th27}

Proof. If $A \in \cd_{G^\s \times L_J^\t}^{cs}(G/U_{P_J})$, then by
Proposition 4.6, $HC_J(A) \in \cd_{G^\s \times L_J^\t}^{cs}(Y_J)$.
Conversely, assume that $HC_J(A) \in \cd_{G^\s \times
L_J^\t}^{cs}(Y_J)$. Then $CH_J \circ HC_J(A) \in \cd_{G^\s \times
L_J^\t}^{cs}(G/U_{P_J})$. Hence by Proposition 4.5, $A \in \cd_{G^\s
\times L_J^\t}^{cs}(G/U_{P_J})$. \qed

\subsection*{4.9} Let $\pi: G/U_{P_J} \rightarrow X_{J, \t}=G/U_{P_J} L_J^\t$
be the quotient map. We call a simple perverse sheaf $A$ in
$\cd_{G^\s}(X_{J, \t})$ a character sheaf on $X_{J, \t}$ if $\tilde
\pi(A) \in \cc_{G^\s \times L_J^\t}(G/U_{P_J})$. They form a set
$\cc_{G^\s}(X_{J, \t})$.

Let $\cd_{G^\s}^{cs}(X_{J, \t})$ be the full subcategory of
$\cd_{G^\s}(X_{J, \t})$ that corresponds to $\cd_{G^\s \times
L_J^\t}^{cs}(G/U_{P_J})$ under $\pi^*$. Then $\cd_{G^\s}^{cs}(X_{J,
\t})$ is the full subcategory of $\cd_{G^\s}(X_{J, \t})$ consisting
of objects whose perverse constituents are contained in
$\cc_{G^\s}(X_{J, \t})$.

\

Now we describe $\cc_{G^\s}(X_{J, \t})$ and $\cd_{G^\s}^{cs}(X_{J,
\t})$ in a slightly different way.

\subsection*{4.10} Define the action of $B$ on $G \times L_J$
by $b \cdot (g, l)=(g b \i, \pi_J(b) l)$. Let $G \times_B L_J$ be
the quotient spaces. Then we may identify $G \times_B L_J$ with
$G/U_{P_J} \times L_J/B_J$ via $(g, l) \mapsto (g l U_{P_J}, l \i
B_J)$. Under this identification, the map $p_J: G \times_B L_J
\rightarrow G/U_{P_J}$ sends $(g, l)$ to $g l U_{P_J}$ and the map
$q_J: G \times_B L_J \rightarrow Y_J$ sends $(g, l)$ to $(g U, l \i
U_J)T$.

Consider the following diagram \[\xymatrix{ G/U_{P_J} \ar[d]_{\pi} &
G \times_B L_J \ar[l]_{p_J} \ar[r]^{q_J} \ar[d]_{i_J} & Y_J \\
X_{J, \t} & G \times_B L_J/L_J^{\t'} \ar[l]_(.6){p'_J} & }\] where
$i_J$ is the projection and $p'_J (g, l)=g l U_{P_J} L_J^\t/U_{P_J}
L_J^{\t}$.

It is easy to see that the square $(p_J, \pi, i_J, p'_J)$ is
Cartesian.

Define (parabolic) character functor $ch_J: \cd_{G^\s \times
L_J^\t}(Y_J) \rightarrow \cd_{G^\s}(X_{J, \t})$ and (parabolic)
Harish-Chandra functor $hc_J: \cd_{G^\s}(X_{J, \t}) \rightarrow
\cd_{G^\s}(Y_J)$ as follows:

For $A \in \cd_{G^\s \times L_J^\t}(Y_J)$, $q_J^* A \in \cd_{G^\s
\times L_J^\t}(G \times_B L_J)$. Let $A'$ be the unique element in
$\cd_{G^\s}(G \times_B L_J/L_J^\t)$ with $\tilde q_J (A)=\tilde i_J
(A')$. Set $ch_J(A)=(p'_J)_! (A')$.

For $B \in \cd_{G^\s}(X_{J, \t})$, set $hc_J(B)=(q_J)_! (p'_J \circ
i_J)^*(B)$.

Using the above diagram, one can easily see that

(a) a simple perverse sheaf in $\cd_{G^\s}(X_{J, \t})$ is a
character sheaf if and only if it is a direct summand of $ch_J(A)$
for some simple perverse sheaf $A \in \cd_{G^\s \times
L_J^\t}^{cs}(Y_J)$;

(b) $A \in \cd_{G^\s}^{cs}(X_{J, \t})$ if and only if $hc_J(A) \in
\bigoplus_{\cl \in \ck(T)} \cd^{\cl}_{G^\s \times L_J^\t}(Y_J)$.

\subsection*{4.11} Notice that each $G^\s \times L_J^\t$-stable subvariety of $Y_J$
that is also stable under the action of $T$ defined in Position 4.6
is of the form $\sqcup_i (G_i/U \times L_i/U_J)/T$, where $G_i$ are
some $G^\s \times B$-orbits on $G$ and $L_i$ are some $L^\t \times
B_J$-orbits on $L_J$. Let $\pi: G \times L_J \rightarrow Y_J$ be the
map defined by $\pi(g, l)=(g U, l \i U_J) T$.  Then one can show
that a simple perverse sheaf in $\cd_{G^\s \times L_J^\t}(Y_J)$ is
contained in $\cd^{\cl}_{G^\s \times L_J^\t}(Y_J)$ if and only if
its image under $\pi^*$ is of the form $A \boxtimes B$, where $A$ is
a simple perverse sheaf in $\cd_{G^\s}(G)$, that is equivariant for
the right $U$-action and has weight $\cl$ for the right $T$-action
and $B$ is a simple perverse sheaf in $\cd_{L_J^\t}(L_J)$, that is
equivariant for the left $U$-action and has weight $\cl \i$ for the
left $T$-action.

Let $pr: G \times L_J/L_J^\t \rightarrow G \times_B L_J/L_J^\t$ be
the projection map. Then a simple perverse sheaf in
$\cd_{G^\s}(X_{J, \t})$ is a character sheaf if and only if it is a
direct summand of $(p'_J)_!(C)$, where $C \in \cd_{G^\s}(G \times_B
L_J/L_J^\t)$ with $\tilde pr (C)=A \boxtimes B'$ for some simple
perverse sheaf $A \in \cd_{G^\s}(G)$, that is equivariant for the
right $U$-action and has weight $\cl$ for the right $T$-action and
some simple perverse sheaf $B \in \cd(L_J/L_J^\t)$, that is
equivariant for the left $U$-action and has weight $\cl \i$ for the
left $T$-action.

\subsection*{4.12} Recall that we have the $G^\s$-stable pieces
decomposition $$X_{J, \t}=\bigsqcup_{w \in \cw(J, \s, \t)}
\bigsqcup_{v_1 \in G^{\s} \setminus G_w/P_K} \bigsqcup_{v_2 \in
L_J^{\t} \setminus L_w/(P_K \cap L_J)} X_{J, \t; v_1, v_2}.$$

Now we define the character sheaves on each piece $X_{J, \t; v_1,
v_2}$. The definition is similar to \cite[4.6]{L3}.

We keep the notation of 3.14. Consider the diagram
\[\xymatrix{ L_K/L_K^{\t'} & G^{\s'}
\times L_K/L_K^{\t'} \ar[l]_(.6){a_1} \ar[r]^{a_2} &
ind^1(L_K/L_K^{\t'}) & ind^1(H L_J^{\t'}/L_J^{\t'}) \ar[l]_{\v} }\]
where $a_1$ and $a_2$ are projections.

Notice that $G^{\s'} \cap P_K=(G^{\s'} \cap U_{P_K}) \rtimes
L_K^{\s'}$, $G^{\s'} \cap U_{P_K}$ is unipotent and $G^{\s'} \cap
U_{P_K}$ acts trivially on $L_K/L_K^{\t'}$. By \cite[A6]{MV},
\begin{align*} \cd_{G^{\s'} \cap
P_K}(L_K/L_K^{\t'})=\cd_{L_K^{\s'}}(L_K/L_K^{\t'}). \tag{a}
\end{align*}

Define the action of $G^{\s'} \times (G^{\s'} \cap P_K)$ on $G^{\s'}
\times L_K/L_K^{\t'}$ by $(g, p) \cdot (g', l)=(g g' p \i, \pi_K(p)
l)$. Then $G^{\s'}$ acts freely on $G^{\s'} \times L_K/L_K^{\t'}$
and $L_K/L_K^{\t'}$ is the quotient space. By \cite[Theorem
2.6.2]{BL}, \begin{align*} \tag{b}  & a_1^*: \cd_{G^{\s'} \cap
P_K}(L_K/L_K^{\t'}) \rightarrow \cd_{G^{\s'} \times (G^{\s'} \cap
P_K)} (G^{\s'} \times L_K/L_K^{\t'}) \\ & \text{ is an equivalence
of categories}. \end{align*}

Similarly, $G^{\s'} \cap P_K$ acts freely on $G^{\s'} \times
L_K/L_K^{\t'}$ and $ind^1(L_K/L_K^{\t'})$ is the quotient space. By
\cite[Theorem 2.6.2]{BL}, \begin{align*} \tag{c}  & a_2^*:
\cd_{G^{\s'}}(Z) \rightarrow \cd_{G^{\s'} \times (G^{\s'} \cap P_K)}
(G^{\s'} \times L_K/L_K^{\t'}) \text{ is an }
\\ & \text{ equivalence of categories}. \end{align*}

\begin{th28} The functors \begin{align*} & \v^*:
\cd_{G^{\s'}}(ind^1(L_K/L_K^{\t'})) \rightarrow
\cd_{G^{\s'}}(ind^1(H L_J^{\t'}/L_J^{\t'})) \\ & \v_!:
\cd^{cs}_{G^{\s'}}(ind^1(H L_J^{\t'}/L_J^{\t'})) \rightarrow
\cd_{G^{\s'}}(ind^1(L_K/L_K^{\t'}))
\end{align*} are equivalences of categories.
\end{th28}

Proof. By Corollary 3.28, $\v$ is an affine space bundle map. Hence
$\v_! \v^*(C)$ is just a shift of $C$ for $C \in
\cd_{G^{\s'}}(ind^1(L_K/L_K^{\t'}))$.

Again by Corollary 3.28, for any $C' \in \cd_{G^{\s'}}(ind^1(H
L_J^{\t'}/L_J^{\t'}))$, $C'$ is constant along each fiber of $\v$.
Hence $\v^* \v_!(C')$ is also a shift of $C'$. The lemma is proved.
\qed

\subsection*{4.14} Combining the above lemma with 4.12 (a), (b) and (c),
we have that the categories $\cd_{L_K^{\s'}}(L_K/L_K^{\t'})$ and
$\cd_{G^{\s'}}(ind^1(H L_J^{\t'}/L_J^{\t'}))$ are naturally
equivalent.

For $X \in \cd^{cs}_{L_K^{\s'}} (L_K/L_K^{\t'})$, let $X'$ be the
unique element in $\cd_{G^{\s'}}(Z)$ such that $a_2^*X'=a_1^* X$.
Set $\tilde X=\v^*(X')$.

As in 4.7, we denote by $\cd^{cs}_{G^{\s'}}(ind^1(H
L_J^{\t'}/L_J^{\t'}))$ the full subcategory of
$\cd_{G^{\s'}}(ind^1(H L_J^{\t'}/L_J^{\t'}))$ whose objects are
$\tilde X$ as above. The simple perverse sheaves that are contained
in $\cd^{cs}_{G^{\s'}}(ind^1(H L_J^{\t'}/L_J^{\t'}))$ are called
character sheaves on $ind^1(H L_J^{\t'}/L_J^{\t'})$.

By Lemma 3.26, $\kappa: ind^1(H L_J^{\t'}/L_J^{\t'}) \xrightarrow
\simeq X_{J, \t; v_1, v_2}$ is an isomorphism. Hence $\kappa^*:
\cd_{G^\s}(X_{J, \t; v_1, v_2}) \rightarrow \cd_{G^{\s'}}(ind^1(H
L_J^{\t'}/L_J^{\t'}))$ is an equivalence of categories. A simple
perverse sheaf $C$ in $\cd_{G^\s}(X_{J, \t; v_1, v_2})$ is called a
character sheaf on $X_{J, \t; v_1, v_2}$ if $\kappa^*(C)$ is a
character sheaf on $ind^1(H L_J^{\t'}/L_J^{\t'})$. We also denote by
$\cd^{cs}_{G^\s}(X_{J, \t; v_1, v_2})$ the subcategory of
$\cd_{G^\s}(X_{J, \t; v_1, v_2})$ that corresponds to
$\cd^{cs}_{G^{\s'}}(ind^1(H L_J^{\t'}/L_J^{\t'}))$ under $\kappa^*$.

\

Now we describe $\cd^{cs}_{G^{\s'}}(ind^1(H L_J^{\t'}/L_J^{\t'}))$
in a different way. This description will be used to prove the main
theorem.

\subsection*{4.15}
Consider the diagram \[\xymatrix{Y & Z_1 \ar[l]_{a} \ar[r]^b & Z_2
\ar[r]^(.3)c & ind^2(H L_J^{\t'}/L_J^{\t'}) \ar[r]^d & ind^1(H
L_J^{\t'}/L_J^{\t'})}\] where $Y$ is the quotient of $ind^2(L_K)/U_K
\times L_J^{\t'} H/U_K$ modulo the diagonal $T$-action on the right,
$Z_1=ind^2(L_K) \times_{B_K} H L_J^{\t'}, Z_2=ind^2(L_K)
\times_{B_K} H L_J^{\t'}/L_J^{\t'}$, $a, b, c$ are analogous to
$q_J, i_J, p_J'$ in 4.10 and $d$ is the projection map.

As in 4.10, we define \begin{align*} & \tilde ch: \cd_{G^{\s'}
\times L_J^{\t'}}(Y) \rightarrow \cd_{G^{\s'}}(ind^2(H
L_J^{\t'}/L_J^{\t'})), \\ & \tilde hc: \cd_{G^{\s'}}(ind^2(H
L_J^{\t'}/L_J^{\t'})) \rightarrow \cd_{G^{\s'} \times L_J^{\t'}}(Y)
\end{align*} as follows:

For $A \in \cd_{G^{\s'} \times L_J^{\t'}}(Y)$, let $A'$ be the
unique element in $\cd_{G^{\s'}}(Z_2)$ with $\tilde b (A')=\tilde a
(A)$. Set $\tilde ch(A)=c_!(A')$.

For $B \in \cd_{G^{\s'}}(ind^2(H L_J^{\t'}/L_J^{\t'}))$, set $\tilde
hc(B)=a_! (c \circ b)^*(B)$.

\begin{th29} Define the action of $U_H$ on $L_J^{\t'} H$ by $g \cdot
g'=g' g \i$. Let $A$ be an $U_H$-equivariant object in
$\cd_{L_J^{\t'}}(L_J^{\t'} H)$ and $\pi: L_J^{\t'} H \rightarrow
L_J^{\t'} H/U_K$ be the projection map. Then $\pi^* \pi_!(A)$ is
also $U_H$-equivariant.
\end{th29}

Proof. Consider the diagram \[\xymatrix{L_J^{\t'} H \times U_K
\times U_H \ar[r]^(.6){m'} \ar[d]_{p_1 \times id} & L_J^{\t'} H
\times
U_K \ar[r]^{p_2} \ar[d]^{p_1} & L_J^{\t'} H \ar[d]^{\pi} \\
L_J^{\t'} H \times U_H \ar[r]^m & L_J^{\t'} H \ar[r]^{\pi} &
L_J^{\t'} H/U_K}\] where $m'(g, l, u)=(g l u \i l \i, l)$, $m(g,
u)=g u \i$, $p_1(g, l)=g l$ and $p_2(g, l)=g$.

It is easy to see that all the squares in the above diagram are
Cartesian.

We have that \begin{align*} m^* \pi^* \pi_!(A) &=m^* (p_1)_! (p_2)^*
(A)=(p_1 \times id)_! (m')^* (p_2)^*(A) \\ &=(p_1 \times id)_!
(m')^*(A \otimes \bar{\mathbb Q}_l^{U_K}). \end{align*}

Consider the diagram \[\xymatrix{L_J^{\t'} H \times U_K \times U_H
\ar[r]^b & L_J^{\t'} H \times U_H \times U_K \ar[r]^{m \times id} &
L_J^{\t'} H \times U_K}\] where $b(g, l, u)=(g, l u l \i, l)$. Then
$m'=b \circ (m \times id)$.

Since $A$ is $U_H$-equivariant, we have that $m^* A \cong A \otimes
\bar{\mathbb Q}_l^{U_H}$. Hence \begin{align*} m^* \pi^* \pi_!(A)
&=(p_1 \times id)_! b^* (m^* A \otimes \bar{\mathbb Q}_l^{U_K})
\cong (p_1 \times id)_! b^* (A \otimes \bar{\mathbb Q}_l^{U_H \times
U_K}) \\ &=(p_1 \times id)_! (A \otimes \bar{\mathbb Q}_l^{U_K
\times U_H})=(p_1)_! (A \otimes \bar{\mathbb Q}_l^{U_K}) \otimes
\bar{\mathbb Q}_l^{U_H} \\ &=(p_1)_! (p_2)^* (A) \otimes
\bar{\mathbb Q}_l^{U_H}=\pi^* \pi_! (A) \otimes \bar{\mathbb
Q}_l^{U_H}.
\end{align*} \qed

\begin{th30} Keep the notation of 4.15. Define $\pi: ind^2(L_K) \times H L_J^{\t'} \rightarrow Y$
by $\pi(z, l)=(z U_K, l \i U_K) T$. The group $U_H$ acts on
$ind^2(L_K) \times H L_J^{\t'}$ on the second factor on the left.
Then for $A \in \cd_{G^{\s'}} (ind^1(H L_J^{\t'}/L_J^{\t'}))$,
$\pi^* \tilde hc(d^* A)$ is $U_H$-equivariant.
\end{th30}

Proof. Consider the following commuting diagram
\[\xymatrix{ind^2(L_K) \times H L_J^{\t'}/L_J^{\t'} \ar[r]^{a_1} \ar[d]_{\pi_1} & ind^2(L_K) \times_{B_K} H L_J^{\t'}/L_J^{\t'}
\ar[r]^(.6){a_2} \ar[d]^{\pi_2} & ind^1(H L_J^{\t'}/L_J^{\t'}) \ar[d]^{\v} \\
ind^2(L_K) \times L_K/L_K^{\t'} \ar[r]^{a_3} & ind^2(L_K)
\times_{B_K} L_K/L_K^{\t'} \ar[r]^(.6){a_4} &
ind^1(L_K/L_K^{\t'})}\] where $a_1, a_3$ are projections, $a_2=d
\circ c$, where $c, d$ are defined in 4.15, $a_4$ is analogous to
$a_2$ and $\pi_1, \pi_2$ are induced from the projection map $H
L_J^{\t'}/L_J^{\t'} \rightarrow L_K/L_K^{\t'}$.

By Lemma 4.13, $A=\v^* (B)$ for some $B \in
\cd_{G^{\s'}}(ind^1(L_K/L_K^{\t'}))$. Then $a_1^* a_2^* (A)=\pi_1^*
(a_4 \circ a_3)^*(B)$. Notice that $\pi_1$ is $U_H$-equivariant,
where $U_H$ acts on $ind^2(L_K) \times L_K/L_K^{\t'}$ trivially.
Hence $a_1^* a_2^*(A)$ is $U_H$-equivariant.

Consider the following commuting diagram \[\xymatrix{ Y'
\ar[d]_{a_5} & ind^2(L_K) \times H L_J^{\t'} \ar[l]_(.7){a'}
\ar[r]^{b'} \ar[d]_{a_6} & ind^2(L_K) \times H L_J^{\t'}/L_J^{\t'}
\ar[d]^{a_1}
\\ Y & ind^2(L_K) \times_{B_K} H L_J^{\t'} \ar[l]_(.7)a \ar[r]^b &
ind^2(L_K) \times_{B_K} H L_J^{\t'}/L_J^{\t'}}\] where
$Y'=ind^2(L_K) \times L_J^{\t'} H/U_K$, $a, b$ are defined in 4.15,
$a'(z, l)=(z, l \i U_K)$, $b'$ is analogous to $b$, $a_5, a_6$ are
projection maps.

It is easy to see that $a_5 \circ a'=\pi$ and the square $(a', a_5,
a_6, a)$ is Cartesian. Now \begin{align*} \pi^* \tilde hc(d^* A)
&=(a')^* a_5^* a_! b^* (a_2)^*(A)=(a')^* (a')_! a_6^* b^* a_2^*(A)
\\ &=(a')^* (a')_! (b')^* a_1^* a_2^*(A). \end{align*}

Since $a_1^* a_2^*(A)$ is $U_H$-equivariant, $(b')^* a_1^* a_2^*(A)$
is also $U_H$-equivariant. Similar to the proof of Lemma 4.16,
$\pi^* \tilde hc(d^* A)$ is $U_H$-equivariant. \qed

\subsection*{4.18} Let $\cd^{cs}_{G^{\s'} \times L_J^{\t'}}(Y)$ be the
full subcategory of $\cd_{G^{\s'} \times L_J^{\t'}}(Y)$ consisting
of elements in $\oplus_{\cl \in \ck(T)} \cd^{\cl}_{G^{\s'} \times
L_J^{\t'}}(Y)$ whose inverse image under the map $\pi$ defined in
the previous Proposition is $U_H$-equivariant.

\begin{th31} A simple perverse sheaf in $\cd_{G^{\s'}}(ind^2(H L_J^{\t'}/L_J^{\t'}))$
is of the form $\tilde d(C)$ for some character sheaf $C$ on
$ind^1(H L_J^{\t'}/L_J^{\t'})$ if and only if it is a direct summand
of $\tilde ch(A)$ for some simple perverse sheaf $A \in
\cd^{cs}_{G^{\s'} \times L_J^{\t'}}(Y)$.
\end{th31}

Proof. Set $X=L_K/L_K^{\t'}$ and $X'=H L_J^{\t'}/L_J^{\t'}$.
Consider the following commuting diagram

\[\xymatrix{L_K \times X \ar[d]_{\pi_1} & L_K \times X'
\ar[l]_{pr_1} \ar[d]_{\pi_2} & Z \ar[l]_(.4){p_1} \ar[r]^(.4){p_2}
\ar[d]^{\pi_3} & ind^2(L_K) \times X' \ar[d]^{\pi_4} \\
L_K \times_{B_K} X \ar[d]_{m} & L_K \times_{B_K} X' \ar[l]_{pr_2}
\ar[d]_{m'} & Z' \ar[l]_(.4){p_3} \ar[r]^(.4){p_4} \ar[d]^{id \times
m'} & ind^2(L_K) \times_{B_K} X' \ar[d]^c \\ X & X' \ar[d]_{pr}
\ar[l]_{pr} & G^{\s'} \times X' \ar[l]_{p_5}
\ar[r]^{p_6} \ar[d]^{id \times pr} & ind^2(X') \ar[d]^{pr_3} \\
& X & G^{\s'} \times X \ar[l]_{p_7} \ar[r]^{p_8} & ind^2(X)}\] where
$Z=G^{\s'} \times L_K \times X'$, $Z'=G^{\s'} \times L_K
\times_{B_K} X'$, $\pi_i, p_i$ are the projection maps, $pr_i$ are
induced from the projection map $pr: X' \rightarrow X$, $m'(l, l')=l
l'$ and $m'$ is analogous to $m$.

It is easy to see that all the squares in the above diagram are
Cartesian.

Let $A$ be a simple perverse sheaf in $\cd^{\cl}_{G^{\s'} \times
L_J^{\t'}}(Y)$. Similar to 4.11, $\tilde \pi(A)=A_1 \boxtimes A_2$,
where $A_1$ is a simple perverse sheaf in
$\cd_{G^{\s'}}(ind^2(L_K))$ that is equivariant for the right
$U_K$-action and has weight $\cl$ for the right $T$-action and $A_2$
is a simple perverse sheaf in $\cd_{L_J^{\t'}}(H L_J^{\t'})$ that is
equivariant for the left $U_K U_H$-action and has weight $\cl \i$
for the left $T$-action.

Consider the diagram \[\xymatrix{L_K & G^{\s'} \times L_K
\ar[l]_{s_1} \ar[r]^{s_2} & ind^2(L_K)}\] where $s_1, s_2$ are
projections. Let $A_1'$ be the unique element in $\cd_{L_K^{\s'}}
(L_K)$ with $\tilde s_1(A_1')=\tilde s_2(A_1)$.

Since $A_2$ is $U_H$-equivariant, there exists a unique simple
perverse sheaf $A'_2$ on $X$ with $\tilde s_3 \tilde pr(A'_2)=A_2$,
where $s_3: H L_J^{\t'} \rightarrow X'$ is the projection map.

Now let $B$ be the simple perverse sheaf in
$\cd_{G^{\s'}}(ind^2(L_K) \times_{B_K} X')$ with $\tilde
\pi_4(B)=A_1 \boxtimes \tilde pr(A'_2)$ and $B'$ be the simple
perverse sheaf in $\cd_{G^{\s'}}(L_K \times_{B_K} X)$ with $\tilde
\pi_1(B')=A'_1 \boxtimes A'_2$. Then $\tilde p_4(B)=\tilde p_3
\tilde pr_2(B')$. We have that \begin{align*} \tilde p_6 \tilde
ch(A) &=\tilde p_6 c_!(B)=(id
\times m')_! \tilde p_4(B)=(id \times m')_! \tilde p_3 \tilde pr_2(B') \\
&=\tilde p_5 m'_! \tilde pr_2(B')=\tilde p_5 \tilde pr \,
m_!(B')=\widetilde{(id \times pr)} \tilde p_7 m_!(B').
\end{align*}

Let $C_1$ be a simple perverse sheaf in $\cd_{G^{\s'}}(ind^2(X'))$
that is a direct summand of $\tilde ch(A)$. Then there exists a
simple perverse sheaf $C'$ in $\cd_{L_K^{\s'}}(X)$ that is a direct
summand of $m_!(B')$ such that $\tilde p_6(C_1)=\widetilde{(id
\times pr)} \tilde p_7(C')$. By 4.11, $C'$ is a character sheaf.

Let $d': ind^2(X) \rightarrow ind^1(X)$ be the projection map. Then
$d' \circ p_8=a_2$, where $a_2$ is defined in 4.12. Let $C_2$ be the
unique element in $\cd_{G^{\s'}}(ind^1(X'))$ with $\tilde
a_2(C_2)=\tilde p_7(C')$. Then $$\widetilde{(id \times pr)} \tilde
p_7(C')=\widetilde{(id \times pr)} \tilde p_8 \tilde d'(C_2)=\tilde
p_6 \tilde d \, \tilde \v(C_2).$$ Therefore $C_1=\tilde d \, \tilde
\v(C_2)$. The ``if'' part is proved.

The ``only if'' part can be proved in the similar way. \qed

\begin{th32} Let $C \in \cd_{G^{\s'}}(ind^1(X'))$. Then $C \in
\cd^{cs}_{G^{\s'}}(ind^1(X'))$ if and only if $\tilde hc
\bigl(\tilde d(C) \bigr) \in \cd^{cs}_{G^{\s'} \times
L_J^{\t'}}(Y)$.
\end{th32}

Proof. Let $Y_0$ be the quotient of $L_K/U_K \times L_J^{\t'} H/U_K$
modulo the diagonal $T$-action on the right. Then $Y=ind^2(Y_0)$,
where $L_K^{\s'}$ acts on $Y_0$ on the first factor. We also have
that $ind^2(L_K) \times_{B_K} H L_J^{\t'}=ind^2(L_K \times_{B_K} H
L_J^{\t'})$. Now consider the following commuting diagram
\[\xymatrix{ Y_0 & L_K \times_{B_K} H L_J^{\t'} \ar[l]_(.6)q \ar[r]^(.6)p &
H L_J^{\t'} \\ G^{\s'} \times Y_0 \ar[u]^{\pi_1} \ar[d]_{\pi_4} &
G^{\s'} \times (L_K \times_{B_K} H L_J^{\t'}) \ar[l]_(.6){id \times
q} \ar[r]^(.6){id \times p} \ar[u]^{\pi_2} \ar[d]_{\pi_5} & G^{\s'}
\times H L_J^{\t'} \ar[u]_{\pi_3} \ar[d]^{\pi_6} \\ Y & ind^2(L_K)
\times_{B_K} H L_J^{\t'} \ar[l]_(.6)a \ar[r]^(.6){c'} \ar[d]_b &
ind^2(H L_J^{\t'}) \ar[d]^{b'} \\ & ind^2(L_K) \times_{B_K} X'
\ar[r]^(.6)c & ind^2(X')}\] where $p, q$ are analogous to $p_J, q_J$
defined in 4.4, $\pi_i$ are the projection maps, $a, b, c$ are
defined in 4.15, $c'$ is induced from $id \times p$, $b': ind^2(H
L_J^{\t'}) \rightarrow ind^2(H L_J^{\t'})/L_J^{\t'}=ind^2(X')$ is
the projection map.

It is easy to see that all the squares in the above diagram are
Cartesian.

Similarly to Proposition 4.5 and 4.6, we can show that

(1) Let $A \in \cd_{L_K^{\s'} \times L_J^{\t'}}(H L_J^{\t'})$, then
some shift of $A$ is a direct summand of $p_! q^* q_! p^*(A)$;

(2) Let $A \in \cd^{\cl}_{L_K^{\s'} \times L_J^{\t'}}(Y_0)$ for some
$\cl \in \ck(T)$, then $q_! p^* p_! q^*(A) \in \sum_{w \in W_K}
\cd^{w^* \cl}_{L_K^{\s'} \times L_J^{\t'}}(Y_0)$.

Now let $C$ be a character sheaf on $ind^1(X')$. Then there exists
$\cl \in \ck(T)$ and a simple perverse sheaf $A \in
\cd^{cs}_{G^{\s'} \times L_J^{\t'}}(Y) \cap \cd^{\cl}_{G^{\s'}
\times L_J^{\t'}}(Y)$ such that $\tilde d(C)$ is a direct summand of
$\tilde ch(A)$. Therefore some shift of $(b')^* \tilde d(C)$ is a
direct summand of $c'_! a^* A$.

It is easy to see that there exists a simple perverse sheaf $A' \in
\cd^{\cl}_{L_K^{\s'} \times L_J^{\t'}}(Y_0)$ with $\pi_1^*
(A')=\pi_4^*(A)$. Then \begin{align*} & \pi_3^* p_! q^*(A')=(id
\times p)_! \pi_2^* q^* (A')=(id \times p)_! (id \times q)^*
\pi_1^*(A') \\ &=(id \times p)_! (id \times q)^* \pi_4^*(A)=(id
\times p)_! \pi_5^* a^*(A)=\pi_6^* c'_! a^*(A). \end{align*}

We can show in the same way that $\pi_1^* q_! p^* p_!
q^*(A')=\pi_4^* a_! (c')^* c'_! a^*(A)$. Hence $\tilde hc(C) \in
\sum_{w \in W_K} \cd^{w^* \cl}_{G^{\s'} \times L_J^{\t'}}(Y)$. By
Proposition 4.17, $\tilde hc(C) \in \cd^{cs}_{G^{\s'} \times
L_J^{\t'}}(Y)$.

On the other hand, for any $C_1 \in \cd_{G^{\s'}}(ind^1(X'))$, there
exists $C'_1 \in \cd_{L_K^{\s'} \times L_J^{\t'}}(H L_J^{\t'})$ with
$\pi_3^*(C'_1)=\pi_6^* (b')^* \tilde{d} (C_1)$. We can also show
that $\pi_3^* p_! q^* q_! p^*(C'_1)=\pi_6^* c'_! a^* a_! (c')^*
(b')^* \tilde d(C_1)=\pi_6^* c'_! a^* \tilde hc \bigl(\tilde
d(C_1)\bigr)$ is a shift of $\pi_6^* (b')^* \tilde ch \circ \tilde
hc \bigl(\tilde d(C_1)\bigr)$.

Hence some shift of $\pi_6^* (b')^* \tilde d (C_1)$ is a direct
summand of $\pi_6^* (b')^* \tilde ch \circ \tilde hc \bigl(\tilde
d(C_1)\bigr)$ and some shift of $\tilde d(C_1)$ is a direct summand
of $\tilde ch \circ \tilde hc \bigl(\tilde d(C_1)\bigr)$. If
moreover, $\tilde hc(C) \in \cd^{cs}_{G^{\s'} \times L_J^{\t'}}(Y)$,
then by Lemma 4.19, $C_1$ is a character sheaf. \qed

\begin{th33} Let $i: X_{J, \t; v_1, v_2} \rightarrow X_{J, \t}$ be
the inclusion map. Then

(1) for any $C \in \cd^{cs}_{G^\s}(X_{J, \t; v_1, v_2})$, $i_! (C),
i_*(C) \in \cd^{cs}_{G^\s}(X_{J, \t})$;

(2) for any $C \in \cd^{cs}(X_{J, \t})$, $i^*(C) \in
\cd^{cs}_{G^\s}(X_{J, \t; v_1, v_2})$.
\end{th33}

Proof. (1) Define the map $m: ind^2(L_K) \rightarrow G$ by $m(g,
l)=g l$ for $g \in G^{\s'}$ and $l \in L_K$. Now consider the
following commuting diagram \[\xymatrix{ X_{J, \t} & ind^2(L_K)
\times_{B_K} X' \ar[l]_f & ind^2(L_K) \times_{B_K} H L_J^{\t'}
\ar[l]_b \ar[d]^a\\ G/U_{P_J} \times L_J/B_J \ar[u]^{i_J \circ p_J}
\ar[d]_{q_J} & Z \ar[l]_(.4){f'} \ar[u]^{\pi_1} \ar[ur]^{\pi_2}
\ar[r]_{\pi_3} \ar[dl]^{\pi_4} & Y \\ Y_J & & }\] where $Z=ind^2(L_K)
\times_{B_K} H L_J^{\t'} \times L_J/B_J$, $f(z, x)=x_1 m(z) x x_2 \i
U_{P_J} L_J^\t$ for $z \in ind^2(L_K)$ and $x \in X'$, $f'(z, h, l
B_J)=(x_1 m(g) h x_2 \i U_{P_J}, l B_J)$ for $z \in ind^2(L_K)$, $h
\in H L_J^{\t'}$ and $l \in L_J$, $\pi_2$ is the projection map,
$\pi_1=b \circ \pi_2$, $\pi_3=a \circ \pi_2$ and $\pi_4=q_J \circ
f'$.

It is easy to see that the square $(\pi_1, f, f', i_J \circ p_J)$ is
Cartesian.

By Lemma 4.19, it suffices to prove that for any $A \in
\cd^{cs}_{G^{\s'} \times L_J^{\t'}}(Y)$ and $B \in
\cd_{G^{\s'}}(ind^2(L_K) \times_{B_K} X')$ with $b^*(B)=a^*(A)$, we
have that $f_!(B), f_*(B) \in \cd^{cs}_{G^\s}(X_{J, \t})$.

Now \begin{align*} hc_J(f_!(B)) &=(q_J)_! (i_J \circ p_J)^*
f_!(B)=(q_J)_! f'_! \pi_1^*(B)=(\pi_4)_! \pi_2^* b^*(B) \\
&=(\pi_4)_! \pi_2^* a^*(A)=(\pi_4)_! (\pi_3)^*(A). \end{align*}

Recall that we have a partition $L_J=\sqcup_{w \in W_J} B_J w B_J$.
Moreover, $B_J w B_J=\sqcup_{i \in \mathbb N \cup \{0\}} L_{w, i}$,
where $$L_{w, i}=\{l \in B_J w B_J; \dim(l \i U_K l \cap U_J)=i\}.$$
It is easy to see that

(a) $L_{w, i}=\varnothing$ for $i \gg 0$;

(b) $L_{w, i}$ is stable under the action of $B_K$ on the left and
the action of $B_J$ on the right;

(c) for any $l \in L_{w, i}$, $l \i U_K l \cap U_J$ is an affine
space of dimension $i$.

Now set $$Z_{w, i}=\{(g, l, l') \in Z; g \in ind^2(L_K), l \in H
L_J^{\t'}, l' \in L_J/B_J, l x_2 \i l' \in L_{w, i}\}.$$

By (b), $Z_{w, i}$ is well-defined. By (a), $Z=\sqcup_{w \in W_J, i
\in \mathbb N \cup \{0\}} Z_{w, i}$ is a finite partition. Hence
$$(\pi_4)_! \pi_3^*(A) \in < (\pi_4 \mid_{Z_{w, i}})_! (\pi_3 \mid_{Z_{w,
i}})^*(A); w \in W_J, i \in \mathbb N \cup \{0\}>.$$

Set \begin{align*} Z'_{w, i}=& \{\bigl((x U_K, y U_K)T, (a U, b
U_J)T \bigr) \in Y \times Y_J; \\ &  U_K m(x) \i x_1 \i a U_J=U_K y
\i x_2 \i b U_J \subset L_{w, i}\}.\end{align*} Define the map
$\pi_{w, i}: Z_{w, i} \rightarrow Z'_{w, i}$ by $$(z, l, l') \mapsto
\bigl( (z U_K, l \i U_K)T, (x_1 m(z) l x_2 \i l' U, l' U_J)T
\bigr)$$ for $z \in ind^2(L_K)$, $l \in H L_J^{\t'}$ and $l' \in
L_J$.

By (c), $\pi_{w, i}$ is an affine space bundle map with fibres
isomorphic to $l \i U_K u \cap U_J$ for $l \in L_{w, i}$.

Let $p_{w, i}: Z'_{w, i} \rightarrow Y$ and $p'_{w, i}: Z'_{w, i}
\rightarrow Y_J$ be the projection maps. Similar to the proof of
Proposition 4.6, $$(\pi_4 \mid_{Z_{w, i}})_! (\pi_3 \mid_{Z_{w,
i}})^*(A)=(p'_{w, i})_! p_{w, i}^*(A) [2 i](i) \in \cd^{w \i
\cl}_{G^\s \times L_J^\t}(Y_J).$$ Hence $hc_J(f_!(B)) \in
\cd^{cs}_{G^\s \times L_J^\t}(Y_J)$. By Proposition 4.8, $f_!(B) \in
\cd^{cs}_{G^\s}(X_{J, \t})$.

We can prove in the same way that $f_* (B) \in \cd^{cs}_{G^\s}(X_{J,
\t})$.

(2) Consider the following commuting diagram
\[\xymatrix{ X_{J, \t} & G \times_B L_J/L_J^\t \ar[l]_{p_J'} & G
\times_B L_J \ar[l]_(.4){i_J} \ar[d]^{q_J} \\ ind^2(L_K)
\times_{B_K} H L_J^{\t'} \ar[u]^{f \circ b} \ar[d]_{a} & Z
\ar[l]_(.4){\pi_2} \ar[u]^{\pi_5} \ar[ur]^{f'} \ar[r]_{\pi_4}
\ar[dl]^{\pi_3} & Y_J
\\ Y & &}\] where $Z, \pi_2, \pi_3, \pi_4$ are defined above, $p'_J$,
$i_J$, $q_J$ are defined in 4.10 and $\pi_5=i_J \circ f'$.

It is easy to see that the square $(\pi_5, p'_J, \pi_2, f \circ b)$
is Cartesian.

Let $A \in \cd^{cs}_{G^{\s} \times L_J^{\t}}(Y_J)$ and $B \in
\cd_{G^\s}(G \times_B L_J/L_J^\t)$ with $i_J^*(B)=q_J^*(A)$. Then
\begin{align*} a_! (f \circ b)^* (p'_J)_! (B)
&=a_! (\pi_2)_! \pi_5^*(B)=(\pi_3)_! (f')^* (i_J)^*(B) \\
&=(\pi_3)_! (f')^* q_J^*(A)=(\pi_3)_! (\pi_4)^*(A). \end{align*}

Similarly to the proof of part (1), we can show that $$a_! (f \circ
b)^* (p'_J)_! (B) \in \bigoplus_{\cl \in \ck(T)} \cd^{\cl}_{G^{\s'}
\times L_J^{\t'}}(Y_J).$$ By Proposition 4.17, $a_! (f \circ b)^*
(p'_J)(B) \in \cd^{cs}_{G^{\s'} \times L_J^{\t'}}(Y)$. Now part (2)
follows from Proposition 4.20. \qed

\begin{th34} Let $C$ be a simple perverse sheaf in $\cd_{G^\s}(X_{J,
\t})$. Then $C$ is a character sheaf if and only if $C$ is the
perverse extension of $C'$, where $C'$ is a character sheaf on some
$G^\s$-stable piece $X_{J, \t; v_1, v_2}$.
\end{th34}

Proof. Since $X_{J, \t}=\sqcup X_{J, \t; v_1, v_2}$, we can find a
$G^\s$-stable piece $X_{J, \t; v_1, v_2}$ such that $\supp(C) \cap
X_{J, \t; v_1, v_2}$ is open dense in $\supp(C)$. Hence $C
\mid_{X_{J, \t; v_1, v_2}}$ is a simple perverse sheaf on $X_{J, \t;
v_1, v_2}$. By part (2) of the above theorem, $C \mid_{X_{J, \t;
v_1, v_2}}$ is a character sheaf on $X_{J, \t; v_1, v_2}$.

On the other hand, if $C'$ is a character sheaf on $X_{J, \t; v_1,
v_2}$. By part (1) of the above theorem, $i_!(C') \in
\cd^{cs}_{G^\s}(X_{J, \t})$, where $i: X_{J, \t; v_1, v_2}
\rightarrow X_{J, \t}$ is the inclusion map. Since the perverse
extension of $C'$ is a quotient of $^pH^0(i_!(C'))$, the perverse
extension of $C'$ is also contained in $\cd^{cs}_{G^\s}(X_{J, \t})$.
Therefore, the perverse extension of $C'$ is a character sheaf on
$X_{J, \t}$. \qed

\section{Lusztig's functors $e^{J'}_J$ and $f^J_{J'}$}

\subsection{} Let $J \subset J' \subset I$ with
$L_J=\t(L_J)$ and $L_{J'}=\t(L_{J'})$.

Define the action of $P_J$ on $G \times (L_{J'} \cap {P_J})/(L_{J'}
\cap {P_J})^\t$ by $p \cdot (g, l)=(g p \i, \pi_{J'}(p) l)$. Let
$Z_{J, J'}$ be the quotient space. Let $\pi: L_{J'} \cap {P_J}
\rightarrow (L_{J'} \cap {P_J})/(L_{J'} \cap U_{P_J}) \cong L_J$ be
the projection. Then $\pi$ induces a morphism $\bar \pi: (L_{J'}
\cap {P_J})/(L_{J'} \cap {P_J})^\t \rightarrow L_J/L_J^\t$. The
morphism $(id, \bar \pi): G \times (L_{J'} \cap {P_J})/(L_{J'} \cap
{P_J})^\t \rightarrow G \times L_J/L_J^\t$ is equivariant under the
${P_J}$-action. Then it induces a morphism $c: Z_{J, J'} \rightarrow
X_{J, \t}$.

The inclusion map $G \times (L_{J'} \cap {P_J})/(L_{J'} \cap
{P_J})^\t \rightarrow G \times L_{J'}/(L_{J'} \cap {P_J})^\t$
induces an isomorphism from $Z_{J, J'}$ to $G \times_{P_{J'}}
L_{J'}/(L_{J'} \cap {P_J})^\t$. Thus the projection map $G
\times_{P_{J'}} L_{J'}/(L_{J'} \cap {P_J})^\t \rightarrow X_{J',
\t}$ induces a map $d: Z_{J, J'} \rightarrow X_{J', \t}$. We will
write $c$ as $c_{J, J'}$ and $d$ as $d_{J, J'}$ if necessary.

Consider the diagram
$$X_{J, \t} \xleftarrow{c} Z_{J, J'} \xrightarrow{d} X_{J', \t}.$$  Define
$$f^J_{J'}: \cd_{G^\s}(X_{J, \t}) \rightarrow \cd_{G^\s}(X_{J', \t}), \quad
e^{J'}_J: \cd_{G^\s}(X_{J', \t}) \rightarrow \cd_{G^\s}(X_{J, \t})$$
by $f^J_{J'}(A)=d_! c^*(A)$, $e^{J'}_J(A')=c_! d^*(A')$.

In the special case where $\mathbf X_{\mathbf J, \s}=Z_{J, 1, G}$
and $\mathbf X_{\mathbf J', \s}=Z_{J', 1, G}$ (see example 3.6), the
functors $e^{\mathbf J'}_{\mathbf J}$ and $f^{\mathbf J}_{\mathbf
J'}$ are just the functors $e^{J'}_J$ and $f^J_{J'}$ defined in
\cite[6.1]{L3}.

\begin{th35} Let $J_1 \subset J_2 \subset J_3 \subset I$ with
$\t(L_{J_i})=L_{J_i}$ for $i=1, 2, 3$. Then $e^{J_2}_{J_1} \circ
e^{J_3}_{J_2}=e^{J_3}_{J_1}$ and $f^{J_2}_{J_3} \circ
f^{J_1}_{J_2}=f^{J_1}_{J_3}$.
\end{th35}

\begin{rmk} This proposition is a generalization of \cite[Lemma 6.2]{L3}. The proof below is similar to the proof in {\it loc.cit}.
\end{rmk}

Proof. Consider the following commuting diagram

\[
\xymatrix{ & & Z_{J_1, J_3} \ar[dl]_{a} \ar[dr]^{b} & &
\\ & Z_{J_1, J_2} \ar[dl]_{c_{J_1, J_2}} \ar[dr]_{d_{J_1, J_2}} & \diamond & Z_{J_2, J_3}
\ar[dl]^{c_{J_2, J_3}} \ar[dr]^{d_{J_2, J_3}} &
\\ X_{J_1, \t} & & X_{J_2, \t} & & X_{J_3, \t}}
\] where $a$ is defined in the similar
way as the map $c$ in 5.1 and $b$ is induced from the projection map
$G \times_{P_{J_3}} L_{J_3}/(L_{J_3} \cap P_{J_1})^\t \rightarrow G
\times_{P_{J_3}} L_{J_3}/(L_{J_3} \cap P_{J_2})^\t$.

Notice that the square $\diamond$ is Cartesian and $c_{J_1, J_2}
\circ a=c_{J_1, J_3}$ and $d_{J_2, J_3} \circ b=d_{J_1, J_3}$. Then
\begin{align*} e^{J_2}_{J_1} \circ e^{J_3}_{J_2} &=(c_{J_1, J_2})_! (d_{J_1,
J_2})^* (c_{J_2, J_3})_! (d_{J_2, J_3})^*=(c_{J_1, J_2})_! a_! b^*
(d_{J_2, J_3})^* \\ &=(c_{J_1, J_3})_! (d_{J_1,
J_3})^*=e^{J_3}_{J_1}, \\ f^{J_2}_{J_3} \circ f^{J_1}_{J_2}
&=(d_{J_2, J_3})_! (c_{J_2, J_3})^* (d_{J_1, J_2})_! (c_{J_1,
J_2})^*=(d_{J_2, J_3})_! b_! a^* (c_{J_1, J_2})^* \\ &=(d_{J_1,
J_3})_! (c_{J_1, J_3})^*=f^{J_1}_{J_3}.
\end{align*} \qed

\begin{th36} We keep the notation of 5.1. Then

(1) If $C \in \cd_{G^\s}^{cs}(X_{J', \t})$, then $e^{J'}_J(C) \in
\cd_{G^\s}^{cs}(X_{J, \t})$.

(2) If $C' \in \cd_{G^\s}^{cs}(X_{J, \t})$, then $f^J_{J'}(C') \in
\cd_{G^\s}^{cs}(X_{J', \t})$.
\end{th36}

\begin{rmk} The special case where $\mathbf X_{\mathbf J, \s}=Z_{J, 1, G}$ and $\mathbf X_{\mathbf J', \s}=Z_{J', 1, G}$
was proved by Lusztig in \cite[6.7(b)]{L3} and \cite[6.4]{L3}.
\end{rmk}

Proof. (1) Define the action of $B$ on $G \times (L_{J'} \cap
{P_J})/(L_{J'} \cap {P_J})^\t$ by $b \cdot (g, l)=(g b \i,
\pi_{J'}(b) l)$. Denote by $G \times_B (L_{J'} \cap {P_J})/(L_{J'}
\cap {P_J})^\t$ the quotient space. We also define $G \times_B
(L_{J'} \cap {P_J})/(L_{J'} \cap U_{P_J})^\t$ and $G \times_B
(L_{J'} \cap {P_J})$ in the similar way.

Consider the following commuting diagram

\[
\xymatrix{ X_{J, \t} & & \ar[ll]_c Z_{J, J'} \ar[rr]^{d} & & X_{J',
\t} \\ G \times_B L_J/L_J^\t \ar[u]^{p'_J} & & Z_1 \ar[ll]_{h}
\ar[u]^{f} \ar[rr]^{k} & & G \times_B L_{J'}/L_{J'}^\t
\ar[u]_{p'_{J'}} \\ G \times_B L_J \ar[u]^{i_J} \ar[drrr]_{q_J} & &
Z_3 \ar[ll]_{r} \ar[u]^{j} & Z_2 \ar[l]_{s}
\ar[r]^{t} \ar[d]^{v} & G \times_B L_{J'} \ar[u]_{i_{J'}} \ar[d]^{q_{J'}} \\
& & & Y_J \ar[r]^{t'} & Y_{J'}}\] where $Z_1=G \times_B (L_{J'} \cap
{P_J})/(L_{J'} \cap {P_J})^\t$, $Z_2=G \times_B (L_{J'} \cap
{P_J})$, $Z_3=G \times_B (L_{J'} \cap {P_J})/(L_{J'} \cap
U_{P_J})^\t$, $h$ (resp. $k$) is defined in the similar way as $c$
(resp. $d$), $f, j, s$ are the projection maps, $t$ is an inclusion
and $t' \bigl((g U, l U_J)T \bigr)=\bigl(g U, l U_{J'} \bigr) T$.

It is easy to see that the squares $(f, c, h, p'_J)$ and $(j, h, r,
i_J)$ are Cartesian. Notice that we may identify $Y_J$ with
$\bigl(G/U \times (L_{J'} \cap P_J)/U_{J'} \bigr)T$ in the natural
way and the map $t'$ is just the inclusion map under this
identification. Thus the square $(t, q_{J'}, v, t')$ is also
Cartesian.

Notice that $s$ is an affine space bundle map. Then
\begin{align*} & hc_J(e^{J'}_J(C))=(q_J)_! (p'_J \circ i_J)^* c_!
d^*(C)=(q_J)_! r_! (f \circ j)^* d^*(C)
\\ &=(q_J \circ r)_! s_! s^* (f \circ j)^* d^*(C)[2 d](d)=v_! (d
\circ f \circ j \circ s)^*(C)[2 d](d) \\ &=v_! (p'_{J'} \circ i_{J'}
\circ
t)^*(C)[2 d](d)=v_! t^* (p'_{J'} \circ i_{J'})^*(C)[2 d](d) \\
&=(t')^* (q_{J'})_! (\bar p_{J'} \circ i_{J'})^*(C)[2 d](d)=(t')^*
hc_{J'}(C)[2 d](d), \end{align*} where $d=\dim \bigl((L_{J'} \cap
U_{P_J})^{\t}\bigr)$.

By Proposition 4.8 $hc_{J'}(C) \in \bigoplus_{\cl \in \ck(T)}
\cd^{\cl}_{G^\s \times L_{J'}^\t}(Y_{J'})$. Therefore,
$hc_J(e^{J'}_J(C))=(t')^* hc_{J'}(C)[2 d](d) \in \bigoplus_{\cl \in
\ck(T)} \cd^{\cl}_{G^\s \times L_J^\t}(Y_J)$. By Proposition 4.8,
$e^{J'}_J(C) \in \cd_{G^\s}^{cs}(X_{J, \t})$.

(2) By 4.10, it suffices to prove that for any simple perverse sheaf
$A \in \cd_{G^\s, \cl}(G \times_B L_J/L_J^\t)$ and $A' \in
\cd^{cs}_{G^\s \times L_J^\t}(Y_J)$ with $i_J^*(A)=q_J^*(A')$, we
have that $f^J_{J'}(p'_J)_!(A) \in \cd_{G^\s}^{cs}(X_{J', \t})$.

We use the above diagram. Then
$$f^J_{J'}(p'_J)_!(A)=d_! c^* (p'_J)_!(A)=
d_! f_! h^*(A)=(p'_{J'})_! k_! h^*(A).$$

Consider the following commuting diagram
\[\xymatrix{ G \times L_J/L_J^\t \ar[d]_{\pi_1} & G \times
(L_{J'} \cap P_J)/(L_{J'} \cap P_J)^\t \ar[l]_(.6){h'}
\ar[r]^(.6){k'} \ar[d]^{\pi_2} & G \times L_{J'}/L_{J'}^\t
\ar[d]^{\pi_3} \\ G \times_B L_J/L_J^\t & G \times_B (L_{J'} \cap
{P_J})/(L_{J'} \cap {P_J})^\t \ar[l]_(.6){h} \ar[r]^(.6){k} & G
\times_B L_{J'}/L_{J'}^\t}\] where $\pi_1, \pi_2, \pi_3$ are
projections, $h'=(id, \bar \pi)$ is defined in 5.1 and $k': G \times
(L_{J'} \cap {P_J})/(L_{J'} \cap {P_J})^\t=G \times (L_{J'} \cap
{P_J}) L_{J'}^\t/L_{J'}^\t \rightarrow G \times L_{J'}/L_{J'}^\t$ is
just the inclusion map.

It is easy to see that all the squares in the diagram are Cartesian.

Then $\pi_3^* k_! h^*(A)=k'_! \pi_2^* h^*(A)=k'_! (h')^*
\pi_1^*(A)$. By 4.11, $\pi_1^*(A)=A_1 \boxtimes A_2$, where $A_1$ is
a simple perverse sheaf in $\cd_{G^\s}(G)$, that is equivariant for
the right $U$-action and has weight $\cl$ for the right $T$-action
and $A_2$ is a simple perverse sheaf in $\cd(L_J/L_J^\t)$, that is
equivariant for the left $U$-action and has weight $\cl \i$ for the
left $T$-action. Thus $k'_! (h')^* \pi_1^*(A)=A_1 \boxtimes A_3$,
where $A_3 \in \cd(L_{J'}/L_{J'}^\t)$ is equivariant for the left
$U$-action, has weight $\cl \i$ for the left $T$-action and is
supported in $(L_{J'} \cap {P_J}) L_{J'}^\t/L_{J'}^\t$. By 4.11,
$f^J_{J'}(p'_J)_!(A) \in \cd_{G^\s}^{cs}(X_{J', \t})$. \qed

\bibliographystyle{amsalpha}

\begin{thebibliography}{BBD}

\bibitem[BBD]{BBD}
A.~A. Beilinson, J.~Bernstein, and P.~Deligne, \emph{Faisceaux
pervers}, Analysis and topology on singular spaces, I (Luminy,
1981), Soc. Math. France, Paris, 1982, pp.~5--171.

\bibitem[BL]{BL} J. Bernstein and V. Lunts, {\em Derived category of equivariant
sheaves}, Lecture Notes in Mathematics 1578, Springer Verlag, 1994.

\bibitem[C]{C}
R. W. Carter, {\em Finite groups of Lie type, Conjugacy classes and
complex characters}, Reprint of the 1985 original. Wiley Classics
Library. A Wiley-Interscience Publication. John Wiley \& Sons, Ltd.,
Chichester (1993).

\bibitem[Gi]{Gi}
V. Ginsburg, {\em Admissible modules on a symmetric space},
Ast\'erisque no.~173-174 (1989), 9, 199--255.

\bibitem[Gr]{Gr}
I.Grojnowski, {\em Character sheaves on symmetric spaces}, Ph.D.
thesis, MIT (1992).

\bibitem[H1]{H1}
X. He, {\em The $G$-stable pieces of the wonderful
compactification}, to appear in Trans. Amer. Math. Soc, math.RT/0412302.

\bibitem[H2]{H2}
X. He, {\em The character sheaves on the group compactification}, to appear in Adv. in Math, math.RT/0508068.

\bibitem[L1]{L1}
G. Lusztig, Character sheaves, I-V, Adv. in Math., 56 (1985), 193-237; 57 (1985), 226-265; 57 (1985), 266-315; 59 (1986), 1-63; 61 (1986), 103-155.

\bibitem[L2]{L2} G. Lusztig, Character sheaves on disconnected groups, I-VII,
Represent. Theory, 7 (2003), 374--403; 8 (2004), 72--124; 8 (2004), 125--144; 8 (2004),
145--178; 8 (2004), 346--376; 8 (2004), 377--413; 9 (2005), 209--266.

\bibitem[L3]{L3}
G. Lusztig, {\em Parabolic character sheaves, I}, Mosc. Math. J.
{\bf 4} (2004), no.~1, 153-179.

\bibitem[L4]{L4}
G. Lusztig, {\em Parabolic character sheaves, II}, Mosc. Math. J.
{\bf 4} (2004), no.~4, 869--896.

\bibitem[MV]{MV} I. Mirkovi\'{c} and K. Vilonen, {\em Characteristic varieties of
character sheaves}, Invent. math. {\bf 93}, 405-418 (1988).

\bibitem[RS]{RS}
R. W. Richardson and T. A. Springer, {\em The Bruhat order on
symmetric varieties}, Geometriae Dedicata {\bf 35} (1990),
389--436.

\bibitem[Sl]{SL}
P. Slodowy, {\em Simple singularities and simple algebraic groups},
Lecture Notes in Mathematics 815, Springer Verlag, 1980.

\bibitem[Sp1]{Sp1}
T. A. Springer, {\em Combinatorics of $B$-orbits in a wonderful
compactification}, Algebraic groups and arithmetic, 99--117, Tata
Inst. Fund. Res., Mumbai, 2004.

\bibitem[Sp2]{Sp2}
T. A. Springer, {\em Remarks on parabolic character sheaves
(preliminary version)}, unpublished.

\bibitem[St]{St}
R. Steinberg, {\it Endomorphisms of linear algebraic groups}, Mem.
Amer. Math. Soc., 80, Amer. Math. Soc., Providence, R.I., 1968.

\end{thebibliography}

\end{document}